\documentclass[10pt]{article}
\usepackage{algpseudocode}
\usepackage{amsmath, amssymb, amsthm}
\usepackage{float, fullpage, graphicx, multirow,parskip, subcaption, setspace}
\usepackage{comment}
\usepackage{url}
\usepackage{enumitem}
\usepackage{tikz}
\usepackage{bm, bbm}
\usetikzlibrary{shapes.geometric, positioning}
\usetikzlibrary{quotes, angles}
\usepackage{rotating}
\usepackage{hyperref}
\usepackage{xr, xr-hyper}
\usepackage{natbib}
\usepackage{placeins}
\usepackage{algorithm2e}
\RestyleAlgo{ruled}

\definecolor{SkyBlue}{RGB}{14, 118, 188}
\definecolor{BrightRed}{RGB}{223,82, 78}

\hypersetup{pdfborder = {0 0 0.5 [3 3]}, colorlinks = true, urlcolor = SkyBlue, linkcolor = BrightRed, citecolor = SkyBlue}

\DeclareMathOperator{\eig}{eig}
\DeclareMathOperator{\tr}{tr}

\DeclareMathOperator*{\argmax}{arg\,max}
\DeclareMathOperator*{\argmin}{arg\,min}

\newtheorem{theorem}{Theorem}

\newtheorem{lemma}{Lemma}
\newtheorem{corollary}{Corollary}

\newtheorem{remark}{Remark}

\newcommand{\skd}[1]{\textcolor{cyan}{\small [skd]: #1}}

\newcommand{\bY}{\bm{Y}}
\newcommand{\by}{\bm{y}}

\newcommand{\bx}{\bm{x}}
\newcommand{\hOmega}{\hat{\Omega}}

\newcommand{\R}{\mathbb{R}}
\newcommand{\N}{\mathcal{N}}
\def\P{\mathbb{P}}
\newcommand{\E}{\mathbb{E}}

\newcommand{\bdelta}{\boldsymbol{\delta}}

\title{Posterior contraction and uncertainty quantification for the multivariate spike-and-slab LASSO}
\author{Yunyi Shen\thanks{Laboratory for Information \& Decision Systems, Massachusetts Institute of Technology.  \url{yshen99@mit.edu}}, and Sameer K. Deshpande\thanks{Department of Statistics, University of Wisconsin--Madison. \url{sameer.deshpande@wisc.edu}} }

\begin{document}
\singlespacing

\maketitle

\begin{abstract}
We study the asymptotic properties of Deshpande et al.\ (2019)'s multivariate spike-and-slab LASSO (mSSL) procedure for simultaneous variable and covariance selection in the sparse multivariate linear regression problem.
In that problem, $q$ correlated responses are regressed onto $p$ covariates and the mSSL works by placing separate spike-and-slab priors on the entries in the matrix of marginal covariate effects and off-diagonal elements in the upper triangle of the residual precision matrix.
Under mild assumptions about these matrices, we establish the posterior contraction rate for the mSSL posterior in the asymptotic regime where both $p$ and $q$ diverge with $n.$
By ``de-biasing'' the corresponding MAP estimates, we obtain confidence intervals for each covariate effect and residual partial correlation.
In extensive simulation studies, these intervals displayed close-to-nominal frequentist coverage in finite sample settings but tended to be substantially longer than those obtained using a version of the Bayesian bootstrap that randomly re-weights the prior.
We further show that the de-biased intervals for individual covariate effects are asymptotically valid.
\end{abstract}

\section{Introduction\label{sec:introduction}}
Suppose we observe $n$ pairs $(\bx_{1}, \by_{1}), \ldots, (\bx_{n}, \by_{n})$ of $p$-dimensional covariate vectors $\bx \in \R^{p}$ and $q$-dimensional response vectors $\by_{i} \in \R^{q}$ from the multi-outcome Gaussian linear regression model
\begin{equation}
\label{eq:general_model}
\by \vert \bx, B,\Omega \sim \N(B^{\top}\bx, \Omega^{-1}).
\end{equation}
Under this model, the $(j,k)$ entry of $B,$ $\beta_{j,k},$ quantifies the \textit{marginal} effect of the $j^{\text{th}}$ covariate on the $k^{\text{th}}$ outcome: for all $x \in \R$, we have 
$$\beta_{j,k} = \E[Y_{k} \vert X_{j} = x+1, X_{-j}] - \E[Y_{k} \vert X_{j} = x, X_{-j}].$$
The matrix $\Omega$ captures the \textit{conditional} dependence relationships between the outcomes after adjusting for the effects of the covariates.
Specifically, if the $(k,k')$ entry of $\Omega,$ $\omega_{k,k'},$ is equal to zero, then we can conclude that after adjusting for the covariates, responses $k$ and $k'$ are conditionally independent.

Fitting the model in Equation~\eqref{eq:general_model} involves estimating $pq + q(q-1)/2$ unknown parameters.
When the total number of parameters exceeds the number of observations $n,$ the parameters $B$ and $\Omega$ are not likelihood-identified.
To make the estimation problem feasible, it is particularly popular and effective to assume that both $B$ and $\Omega$ are sparse.
Arguing from a Bayesian perspective, \citet{Deshpande2019} introduced the multivariate spike-and-slab LASSO (mSSL) that put spike-and-slab LASSO priors \citep{RockovaGeorge2018_ssl} on all entries in $B$ and off diagonal entries in $\Omega.$ 
They then proposed an Expectation Conditional Maximization \citep{MengRubin1993_ecm} algorithm to approximate the \textit{maximum a posteriori} (MAP) estimates of $B$ and $\Omega.$
That algorithm iteratively solved a sequence of penalized maximum likelihood problems in which each parameter $\beta_{j,k}$ and $\omega_{k,k'}$ is individually penalized.
The underlying Bayesian hierarchical model enables \textit{adaptive} penalty mixing so that smaller parameter estimates are shrunk to zero more aggressively than larger parameter estimates.
Thanks to its automatic, adaptive penalization, the mSSL tends to return estimates that are less biased than those returned by fixed penalty alternatives.
See \citet{GeorgeRockova2020_penaltymix} and \citet{Bai2020_ssl_review} for discussion of this general phenomenon. 

Despite the mSSL's excellent empirical performance in finite samples, \citet{Deshpande2019} did not examine its asymptotic properties nor did they attempt to quantify the uncertainty in their final parameter estimates. 
In this paper, we show that with some minor modifications, the mSSL posterior concentrates.
That is, assuming that the observed data were truly generated from the model in Equation~\eqref{eq:general_model} with a sparse $B_{0}$ and $\Omega_{0}$, the mSSL posterior distributions of $B$ and $\Omega$ collapse to point-masses at $B_{0}$ and $\Omega_{0}.$
Our main results, Theorem~\ref{thm:posterior_contraction_mssl} and Corollary~\ref{cor:B_Psi_recovery} establish the rate of such posterior collapse (i.e.\ the posterior concentration rate) for the parameters $B$ and $\Omega$ as well as evaluations of the regression function $XB$. 
We leverage these results to derive asymptotically-valid confidence intervals for $B$ and $\Omega$ using a de-biasing technique.

The remainder of this paper is organized as follows.
After a brief review of the mSSL in Section~\ref{sec:background}, we present our main theoretical results in Section~\ref{sec:theory}.
Then, in Section~\ref{sec:comparisons}, we discuss how our results generalize many known contraction rates for the sparse single-outcome regression model and the sparse Gaussian graphical model.
In Section~\ref{sec:debias}, we derive asymptotically-valid confidence intervals for individual parameters $\beta_{j,k}$ and $\omega_{k,k'}$ using a de-biasing argument.
We assess the finite-sample coverage of these intervals using a simulation study, which we describe in Section~\ref{sec:simulation}.
We conclude with a discussion of limitations and potential future directions in Section~\ref{sec:discussion}.


\section{Background\label{sec:background}}

\subsection{The mSSL prior}

\citet{Deshpande2019} specified independent \textit{continuous} spike-and-slab priors on the entries of $B$ and $\Omega.$
Specifically, they introduced an additional parameter $\theta \in [0,1]$, and modeled each entry of $B$ as conditionally independent draws from a $\text{Laplace}(\lambda_{0})$ distribution (with probability $1-\theta$) or a $\text{Laplace}(\lambda_{1})$ distribution (with probability $\theta$), where $0 < \lambda_{1} \ll \lambda_{0}$ are fixed positive constants.
Since $\lambda_{0} \gg \lambda_{1},$ the $\text{Laplace}(\lambda_{0})$ distribution (the ``spike'') is much more tightly concentrated around zero than the $\text{Laplace}(\lambda_{1})$ distribution (the ``slab'').
The parameter $\theta$ controls the overall proportion of entries of $B$ that are drawn from the slab and can be interpreted as ``large'' or ``significant.''
They specified a $\text{Beta}(a_{\theta}, b_{\theta})$ for $\theta$ to model their uncertainty about that proportion. 

They similarly introduced a further parameter $\eta \in [0,1],$ which governs the proportion of off-diagonal entries in $\Omega$ drawn from a $\text{Laplace}(\xi_{0})$ spike or $\text{Laplace}(\xi_{1})$ slab with fixed constants $0 < \xi_{1} \ll \xi_{0}.$
They modeled each diagonal entry as being drawn from the $\text{Laplace}(\xi_{1})$ slab.
They additionally constrained the prior to the positive definite cone.

The resulting mSSL prior density can be decomposed as $\pi(B,\theta, \Omega, \eta) = \pi(B \vert \theta)\pi(\Omega \vert \eta)\pi(\theta)\pi(\eta)$ where
\begin{align}
\pi(B \vert \theta) &\propto \prod_{j = 1}^{p}{\prod_{k = 1}^{q}{\left(\theta\lambda_{1}e^{-\lambda_{1}\lvert \beta_{j,k}\rvert} + (1-\theta)\lambda_{0}e^{-\lambda_{0}\lvert \beta_{j,k} \rvert}\right)}} \label{eq:B_prior_density} \\
\pi(\Omega \vert \eta) &\propto \mathbbm{1}(\Omega \succ 0) \times \prod_{k = 1}^{q}{\left[e^{-\xi_{1}\omega_{k,k}} \times \prod_{k' > k}{\left[ \eta\xi_{1}e^{-\xi_{1}\lvert \omega_{k,k'} \rvert} + (1-\eta)\xi_{0}e^{-\xi_{0}\lvert \omega_{k,k'} \rvert} \right]} \right]} \label{eq:Omega_prior_density} \\
\pi(\theta)\pi(\eta) &\propto \theta^{a_{\theta} - 1}(1 - \theta)^{b_{\theta} - 1}\eta^{a_{\eta} -1}(1 - \eta)^{b_{\eta} - 1}. \label{eq:theta_eta_prior_density}
\end{align}

\subsection{MAP estimation via Expectation Conditional Maximization}

Letting $X$ be the $n \times p$ design matrix whose $i^{\text{th}}$ row is $\bx_{i}^{\top},$ the model in Equation~\eqref{eq:general_model} yields the likelihood function 
\begin{equation}
\label{eq:likelihood}
p(\bY \vert X,B, \Omega)\propto \lvert\Omega\rvert^{\frac{n}{2}}\exp\left\{-\frac{1}{2}\tr\left[(Y-XB)^\top\Omega(Y-XB)\right]\right\}.
\end{equation}

The likelihood function combined with the prior densities in Equations~\eqref{eq:B_prior_density}--\eqref{eq:theta_eta_prior_density} together determine the joint posterior distribution of $(B,\theta,\Omega,\eta) \vert \bY.$
Unfortunately, the posterior is analytically intractable, high-dimensional, and highly-multimodal, rendering stochastic search techniques like Markov chain Monte Carlo (MCMC) computationally infeasible.
Further, thanks to the non-concavity of the log prior density, direct optimization is generally difficult.
Instead, \citet{Deshpande2019} introduced an Expectation Conditional Maximization \citep{MengRubin1993_ecm} algorithm to target the \textit{maximum a posteriori} (MAP) estimate of $(B, \theta, \Omega, \eta).$

The algorithm works by introducing additional spike-and-slab indicators and then iteratively optimizing the log posterior density after marginalizing out these indicators. 
Formally, it introduces indicators $\bdelta = \{\delta_{k,k'}: 1 \leq k < k' \leq q\},$ one for each off-diagonal element in $\Omega$'s upper-triangle, with each $\delta_{k,k'} \sim \text{Bernoulli}(\eta)$ independently, and 
$$
\pi(\Omega \vert \bdelta) \propto \mathbbm{1}(\Omega \succ 0) \times \prod_{k = 1}^{q}{\left[e^{-\xi_{1}\omega_{k,k}} \times \prod_{k' > k}{\left(\xi_{1}e^{-\xi_{1}\lvert\omega_{k,k'}\rvert}\right)^{\delta_{k,k'}}\left(\xi_{0}e^{-\xi_{0}\lvert \omega_{k,k'}\rvert}\right)^{1 - \delta_{k,k'}}}\right]}.
$$

The algorithm then iterates between an ``E step'' in which it computes a new surrogate objective and two conditional maximization steps (``CM steps'') in which it maximizes that surrogate.
Letting $(B^{(t-1)}, \theta^{(t-1)}, \Omega^{(t-1)}, \eta^{(t-1)})$ be the parameter estimates at the start of the $t^{\text{th}}$ iteration, the E step computes
$$
F^{(t)}(B, \theta, \Omega, \eta) = \E_{\bdelta \vert \cdot}\left[\log \pi(B, \theta, \Omega, \eta, \bdelta \vert \bY) \vert B = B^{(t-1)}, \theta = \theta^{(t-1)}, \Omega = \Omega^{(t-1)}, \eta = \eta^{(t-1)}, \bY\right],
$$
where the expectation is taken with respect to the conditional posterior distribution of the spike-and-slab indicators.
Then, in the CM steps, we sequentially solve the following problems
\begin{align}
(B^{(t+1)}, \theta^{(t+1)}) &= \argmax_{B,\theta} F^{(t)}(B, \theta, \Omega^{(t-1)}, \eta^{(t-1)}) \label{eq:CM_step1} \\
(\Omega^{(t+1)}, \eta^{(t+1)}) &= \argmax_{\Omega,\eta} F^{(t)}(B^{(t)}, \theta^{(t)}, \Omega, \eta) \label{eq:CM_step2}
\end{align}

It turns out that one can compute $B^{(t+1)}$ using a cyclical coordinate descent algorithm that soft-thresholds running estimates of each $\beta_{j,k}$ at $\lambda^{\star}_{j,k}$ 
The entry-specific threshold $\lambda^{\star}_{j,k}$ lies between the smaller slab penalty $\lambda_{1}$ and the larger spike penalty $\lambda_{0}$.
Indeed, $\lambda_{j,k}^{\star} = \lambda_{1}p^{\star}_{j,k} + \lambda_{0}(1 - p^{\star}_{j,k}),$ where $p^{\star}_{j,k}$ is the conditional posterior probability that the previous estimate $\beta^{(t-1)}_{j,k}$ was drawn from the slab distribution.
In this way, the algorithm applies a larger (resp.\ smaller) penalty to the parameter $\beta_{j,k}$ when the previous running estimate is small (resp.\ large) in absolute value and less (resp.\ more) likely to be drawn from the slab distribution.
Computing $\Omega^{(t+1)}$ reduces to solving a graphical LASSO \citep{Friedman2008} problem with individualized penalties on the parameters $\omega_{k,k'}$'s.
These individual penalties similarly range from $\xi_{1}$ to $\xi_{0},$ depending on how likely it was to draw $\omega^{(t-1)}_{k,k'}$ from the slab.
We refer to \citet[\S2]{Deshpande2019} for a detailed derivation and description of the mSSL algorithm.

\citet{Deshpande2019} recommended $a_{\theta} = 1, b_{\theta} = pq, a_{\eta} = 1$ and $b_{\eta} = q,$ choices that strongly encourage sparsity in $B$ and $\Omega.$
They also recommended setting $\lambda_{1} = 1$ and $\xi_{1} = n/100.$
Rather than fixing a single $\lambda_{0}$ and single $\xi_{0}$ value, they specified grids of increasing spike parameters and ran their ECM algorithm for each combination with warm-starts.
This strategy, which \citet{RockovaGeorge2018_ssl} termed ``dynamic posterior exploration,'' essentially enabled sequential maximization of several posterior distributions, one for each combination of $(\lambda_{0}, \xi_{0})$ value.
By increasing $\lambda_{0}$ and $\xi_{0}$ --- thereby making the spike distributions ``spikier'' --- their computational approach propagates estimates of $B$ and $\Omega$ through a series of increasingly strict filters that remove negligible covariate effects and partial covariances.
They specifically recommended $\lambda_{0}$ range from 1 to $n$ and $\xi_{0}$ to range from $n/10$ to $1.$
Though these values are somewhat arbitrary, they appeared to work quite well in practice.

\subsection{Spike-and-slab asymptotics \& uncertainty quantification}

The mSSL algorithm generalizes the SSL algorithm of \citet{RockovaGeorge2018_ssl}, which was proposed for single-outcome (i.e., $q = 1$) high-dimensional regression.
In that work, the authors showed that when the residual variance is known, the posterior distribution of the regression coefficients contracts at a nearly the minimax-optimal rate as $n \rightarrow \infty.$ 
In fact, even when the residual variance is unknown, the posterior distribution of the regression coefficients contracts at a near-minimax rate \cite[Theorem 2 and Remark 2]{Bai2020_groupSSL}.
For Gaussian graphical model estimation, \citet{Gan2019_unequal} showed that the MAP estimator corresponding to placing spike-and-slab LASSO priors on the off-diagonal elements of a precision matrix is consistent. 
\citet{Ning2020} showed that the joint posterior distribution of $(B,\Omega)$ in Equation~\eqref{eq:likelihood} contracts to the true value when using a group spike-and-slab prior with Laplace slab and point mass spike on $B$ and a carefully selected prior on the eigendecomposition of $\Omega^{-1}.$ 
To our knowledge, however, asymptotic results for mSSL have not yet been established.

Posterior contraction results notwithstanding, quantifying the finite sample uncertainty in spike-and-slab posteriors is generally very difficult.
Approaches based on MCMC tend to navigate the high-dimensional, multi-modal posterior landscapes extremely slowly, effectively rendering them computationally prohibitive.
In recent years, several authors have suggested variants of the weighted likelihood bootstrap \citep{Newton1994} for simulating approximate posterior samples \citep[see, e.g.,][]{Newton2021, Nie2020, Menacher2023_bless}.
The weighted likelihood bootstrap works by repeatedly maximizing a weighted version of the log-likelihood function, where the weights for each observation are drawn from a $\text{Gamma}(1,1)$ distribution.
Later, \citet{Newton2021} modified that procedure by re-weighting both the individual likelihood terms and the log-prior.
When run with sparsity-inducing priors like the SSL, repeatedly solving the weighted MAP problem produces exactly sparse samples, which are not representative of the absolutely continuous posterior distribution.
To counter-act this behavior, working in the context of single-outcome regression with the SSL, \citet{Nie2020} additionally randomly re-centered the prior density.
They further showed that the distribution of their re-weighted and re-centered Bayesian bootstrap samples contracts with the same rate as the true posterior towards the true parameters.

In the context of the mSSL, we consider a two-step approach.
First, we run the dynamic posterior exploration strategy using \citet{Deshpande2019}'s recommend settings to obtain estimates $\hat{B}, \hat{\Omega}, \hat{\theta},$ and $\hat{\eta}.$
Then, we repeatedly solve the randomized optimization problem
\begin{equation}
\label{eq:BB_main_equation}
\argmax_{B,\Omega}\left\{\sum_{i = 1}^{n}{w_{i}\log p(\by_{i} \vert \bx_{i}, B,\Omega)} + w_{0}\left[\log\left(\pi(B - B' \vert \hat{\theta})\pi(\Omega - \Omega' \vert \hat{\eta})\right)\right]\right\}.
\end{equation}
where ($w_{1}, \ldots, w_{n})$ are independent $\text{Gamma}(1,1)$ random variables; $B'$ and $\Omega'$ are (possibly random) re-centering matrices; and $\pi(\cdot \vert \hat{\theta})$ and $\pi(\cdot \vert \hat{\eta})$ are the conditional prior densities from Equations~\eqref{eq:B_prior_density} and~\eqref{eq:Omega_prior_density} evaluated at $\theta = \hat{\theta}, \lambda_{1} = 1, \lambda_{0} = n,$ and $\eta = \hat{\eta}, \xi_{1} = n/100, \xi_{0} = 1.$
By setting $w_{0} = 0, B' = \mathbf{0}_{p \times q}$ and $\Omega' = \mathbf{0}_{q \times q},$ we obtain the mSSL analog of \citet{Newton1994}'s weighted likelihood bootstrap.
Similarly, by drawing $w_{0} \sim \text{Gamma}(1,1)$ independently of the the other $w_{i}$'s and setting $B'$ and $\Omega'$ to zero, we obtain an analog to \citet{Newton2021}'s weighted Bayesian bootstrap.
Finally by drawing $w_{0} \sim \text{Gamma}(1,1)$ independently of the other $w_{i}$'s and drawing $B'$ and $\Omega'$ from the appropriate slab distributions, we get the mSSL analog of \citet{Nie2020}'s procedure.
Although these approaches are conceptually attractive --- they can produce independent samples of $(B,\Omega)$ in an embarrassingly parallel fashion ---  we are unaware of any general guarantee the bootstrap samples closely approximate independent samples from the actual posterior.
In fact, we are also unaware of any general result guaranteeing that Bayesian bootstrapping produces uncertainty intervals with nominal frequentist coverage in finite samples.

Instead of relying on a bootstrap-like procedure, \citet{Bai2020_groupSSL} leveraged their posterior contraction results to derive asymptotically valid confidence intervals about high-dimensional regression coefficients using the ``de-biasing'' approach of \citet{van2014asymptotically,javanmard2018debiasing}. 
This approach works by (i) forming a new estimator that removes estimation bias introduced by the penalty and (ii) basing inference on the new estimator's asymptotic sampling distribution.

In summary, asymptotic results for the mSSL have not yet been established nor have any authors attempted to quantify the uncertainty around the mSSL estimator, motivating our work.

\section{Asymptotic theory of mSSL\label{sec:theory}}
To analyze the mSSL theoretically we require some additional assumptions about $B$ and $\Omega$ and must slightly modify the prior.
First, we assume that there is a sparse $B_{0}$ and a sparse $\Omega_{0}$ such that multi-outcome Gaussian linear regression model in Equation~\eqref{eq:general_model} is correctly specified.
Let $s_{0}^{B}$ and $s_{0}^{\Omega}$ respectively be the number of non-zero free parameters in $B_{0}$ and $\Omega_{0}.$
Note that $s_{0}^{B} \leq pq$ and $s_{0}^{\Omega} \leq Q,$ where $Q = q(q-1)/2.$
In what follows, we denote the contribution of the $i^{\text{th}}$ observation to the likelihood function evaluated at a generic $B$ and $\Omega$ (resp. the true $B_{0}$ and $\Omega_{0})$ by $f_{i}$ (resp. $f_{0,i}$).
We additionally denote the full data likelihoods evaluated the same parameters by $f_{0} = \prod_{i=1}^{n}{f_{0,i}}$ and $f = \prod_{i=1}^{n}{f_{i}}.$
For sequences of real numbers $\{a_{n}\}$ and $\{b_{n}\},$ we write $a_{n} \lesssim b_{n}$ to mean that there is a constant $C$ such that for $n$ large enough, $a_{n} \leq C b_{n}.$
Our main assumptions are:
\begin{description}
\item[A1]{$\Omega_{0}$ has bounded operator norm and eigenvalues; that is, 
$$
\Omega_{0} \in \mathcal{H}_{0} = \{\Omega: \text{all of $\Omega$'s eigenvalues lie in $[b_{1}^{-1}, b_{2}^{-1}]$}\},
$$ where $b_{1} > b_{2}$ are fixed, positive constants not depending on $n.$}
\item[A2]{$B_{0}$ has bounded entries; that is, $B_{0} \in \mathcal{B}_{0} = \{B: \lVert B \rVert_{\infty} < a_{1}\}$ where $a_{1} > 0$ is a fixed constant not depending on $n.$}
\item[A3]{Dimensionality:  $\log(n) \lesssim \log(q), \log(n) \lesssim \log(p)$ and $\max\{q,s_0^\Omega,s_0^B\}\log(\max\{p,q\})/n\to 0.$}
\item[A4]{Prior tuning for $B:$ there are constants $a' > 0$ and $b' > 1/2$ that do not depend on $n$ such that (i) $(1-\theta)/\theta\sim (pq)^{2+a'}$; (ii) $\lambda_0\sim \max\{n,pq\}^{2+b'}$; and (iii) $\lambda_1\asymp  1/n.$}
\item[A5]{Prior tuning for $\Omega:$ there are constants $a,b > 0$ not depending on $n$ such that (i) $(1-\eta)/\eta\sim Q^{2+a}$; (ii) $\xi_0\sim \max\{Q,n\}^{4+b}$; and (iii) $\xi_1\asymp 1/\max\{Q,n\}$.}
\end{description}

Before proceeding, we highlight the difference between the assumptions here and the model reviewed in Section~\ref{sec:background}.
First, instead of truncating the the prior for $\Omega$ to the positive definite cone, we need to assume that the smallest eigenvalue of $\Omega$ is bounded away from 0. 
This assumption ensures that the residual matrix $Y - XB$ is not too poorly conditioned, essentially permitting us to ignore the case of substantially overfitting the outcomes $Y$ in our theoretical analysis.
Additionally, we assume that $\theta$ and $\eta$ are fixed and known (Assumptions A4 and A5), rather than being modeled as in the Section~\ref{sec:background}.
Such assumptions mirror those made by \citet{RockovaGeorge2018_ssl}.


\subsection{Posterior contraction of the mSSL}

Our first results (Theorem~\ref{thm:log-affinity}) shows that the mSSL posterior contracts in log-affinity.
\begin{theorem}[Contraction in log affinity]
    \label{thm:log-affinity}
    Under Assumptions A1--A5, there is some constant $M>0$ that does not depend on $n$ such that
    \begin{align}
        \sup_{B\in\mathcal{B}_0,\Omega\in\mathcal{H}_0}\E_0 \Pi\left(B,\Omega:\frac{1}{n}\sum_{i = 1}^{n} \rho(f_i,f_{0i})\ge M\epsilon_n^2|Y_1,\dots,Y_n\right) &\longrightarrow 0
    \end{align}
    where $\epsilon_n=\sqrt{\max\{q,s_0^\Omega,s_0^B\}\log(\max\{p,q\})/n};$ $s_{0}^{B}$ and $s_{0}^{\Omega}$ are the numbers of non-zero entries in $B$ and off-diagonal entries in $\Omega;$ and the log-affinity is defined as
    \begin{align*}
    \frac{1}{n}\sum_{i}^{n} \rho(f_i,f_{0,i}):=&-\log\left(\frac{|\Omega^{-1}|^{1/4}|\Omega_0^{-1}|^{1/4}}{|(\Omega^{-1}+\Omega_0^{-1})/2|^{1/2}}\right)+\frac{1}{8n}\sum_{i = 1}^{n} \bx_i^{\top}(B-B_0)\left(\frac{\Omega^{-1}+\Omega_0^{-1}}{2}\right)^{-1}(B-B_0)^{\top}\bx_{i}.
\end{align*}
\end{theorem}

The proof of Theorem~\ref{thm:log-affinity} follows the standard recipe for deriving posterior concentration rates based on independent but non-identically distributed observations, \textit{viz.,} \citet[Theorem 8.23]{ghosal2017fundamentals}.
We divided our proof into several intermediate lemmas whose proofs, which involve mostly straightforward calculations, we defer to the Supplemental Materials.
Before stating these lemmas, we need to define a few more quantities.
First, let $K$ and $V$ respectively be the Kullback-Leibler divergence and variance between a multivariate regression model with parameters $(B_{0},\Omega_{0})$ and parameters $(B,\Omega)$:
{ \small
\begin{align}
n^{-1}K(f_{0},f) &= \frac{1}{2}\left(\log\left(\frac{|\Omega_0|}{|\Omega|}\right)-q+\tr(\Omega_0^{-1}\Omega)+\frac{1}{n}\tr\left\{ (B_0-B)^{\top}X^{\top}X(B_0-B)\Omega \right\}\right) \label{eq:kldiv} \\
n^{-1}V(f_{0}, f) &= \frac{1}{2}\left(\tr((\Omega_0^{-1}\Omega)^2)-2\tr (\Omega_0^{-1}\Omega) + q \right)+\frac{1}{n}\tr\left((B-B_0)^{\top}X^{\top}X(B-B_0)\Omega\Omega_0^{-1}\Omega\right)\label{eq:klvar}.
\end{align}
}
Given the values of $\theta$ and $\eta,$ let $\delta_{\beta}$ and $\delta_{\omega}$ denote the points at which the spike density and slab densities intersect:
\begin{align*}
\delta_{\beta} &= \frac{1}{\lambda_0-\lambda_1}\log\left[\frac{1-\theta}{\theta}\times \frac{\lambda_0}{\lambda_1}\right] &\text{and} & & \delta_{\omega}=\frac{1}{\xi_0-\xi_1}\log\left[\frac{1-\eta}{\eta}\times\frac{\xi_0}{\xi_1}\right].
\end{align*}
We define the \textit{effective dimensions} of the matrices $B$ and $\Omega$ to be the number of free parameters that exceed $\delta_{\beta}$ and $\delta_{\omega}$:
\begin{align*}
\nu_{\beta}(B) &= \sum_{j,k}{\mathbbm{1}(\lvert \beta_{j,k} \rvert > \delta_{\beta})} &\text{and} & & \nu_{\omega}(\Omega) &= \sum_{k < k'}{\mathbbm{1}(\lvert\omega_{k,k'}\rvert > \delta_{\omega})}.
\end{align*}
 
Lemma~\ref{lemma:KL_m} shows that the modified mSSL prior places enough probability mass in small neighborhoods around every possible choice of $(B_{0},\Omega_{0})$.
\begin{lemma}[KL condition]
\label{lemma:KL_m}
Let $\epsilon_n=\sqrt{\max\{q,s_0^\Omega,s_0^B\}\log(\max\{p,q\})/n}.$ Then for all true parameters $(B_{0},\Omega_{0})$ we have
$$
-\log\Pi\left[(B,\Omega):K(f_0,f)\le n\epsilon_n^2, V(f_0,f)\le n\epsilon_n^2\right]\le C_1n\epsilon_n^2.
$$
Further, consider the event $E_{n} = \left\{\bY: \iint f(\bY)/f_{0}(\bY) d\Pi(B)d\Pi(\Omega) \geq e^{-C_{1}n\epsilon_{n}^{2}}\right\}.$
Then for all $(B_{0}, \Omega_{0}),$ we have $\P_{0}(E_{n}^{c}) \to 0$ as $n \to \infty.$
\end{lemma}

\begin{proof}
Observe that both KL divergence~\eqref{eq:kldiv} and KL variance~\eqref{eq:klvar} has term that only depends on $\Omega$ and term depends on both $B$ and $\Omega$. We separate the two type of terms and bound them separately.

To simplify the notation, we denote $(B-B_0)=\Delta_B$ and $\Omega-\Omega_0=\Delta_\Omega$. 
Consider the events $\mathcal{A}_{1}$ and $\mathcal{A}_{2}$ defined as
\begin{equation}
    \label{eqn:mssl_A1_due_to_Omega}
    \mathcal{A}_1=\left\{ \Omega: \tr((\Omega_0^{-1}\Omega)^2)-2\tr (\Omega_0^{-1}\Omega) + q\le \epsilon_n^2, \log\left(\frac{|\Omega_0|}{|\Omega|}\right)-q+\tr(\Omega_0^{-1}\Omega)\le \epsilon_n^2 \right\}
\end{equation}
and
\begin{equation}
    \label{eqn:mssl_A2_condition_on_A1}
        \mathcal{A}_2=\left\{ (\Omega, B): \tr\left( (B_0-B)^{\top}X^{\top}X(B_0-B)\Omega \right)\le n\epsilon^2, \tr\left((B-B_0)^{\top}X^{\top}X(B-B_0)\Omega\Omega_0^{-1}\Omega\right)  \le n\epsilon_n^2/2  \right\}
\end{equation}

By construction the event $\mathcal{A}_1$ involves only $\Delta_\Omega$ and the event $\mathcal{A}_2$ involves both $\Delta_\Omega$ and $\Delta_B.$  
Further, the event $(\mathcal{A}_1\cap \{\Omega\succ \tau I\})\cap \mathcal{A}_2$ is a subset of the event  $\{ K/n\le \epsilon_n^2, V/n\le\epsilon_n^2 \}$. To show the lemma, it suffices to lower bound $\log\Pi(\mathcal{A}_1)+\log\Pi(\mathcal{A}_1|\mathcal{A}_2)$ with the stated lower bound. We separately bound the prior probabilities $\Pi(\mathcal{A}_1)$ and $\Pi(\mathcal{A}_1|\mathcal{A}_2)$.

We start with lower bound the prior mass put on $\mathcal{A}_1$. We consider a smaller set, the vectorized L1 ball $\mathcal{A}_1^\star$:
\begin{align*}
    \mathcal{A}_1^\star=\{2\sum_{k>k'} |\omega_{0,k,k'}-\omega_{k,k'}|+\sum_k |\omega_{0,k,k}-\omega_{k,k}|\le \frac{\epsilon_n}{b_2}\}
\end{align*}
where $b_{2} > 0$ is the constant in the assumption of $\Omega_0$'s spectra.

Since the Frobenius norm is bounded by the vectorized L1 norm, we conclude that 
$$
\mathcal{A}_1^{\star} \subset \left\{ \lVert \Omega_{0} - \Omega \rVert_{F} \leq \frac{\epsilon_{n}}{b_{2}}\right\}.
$$
We now show that $\left\{ \lVert \Omega_{0} - \Omega \rVert_{F} \leq \frac{\epsilon_{n}}{b_{2}} \right\}\subset \mathcal{A}_1$.

Per \citet{Ning2020}'s Lemma 5.1, using our assumption that $\Omega_0$'s eigenvalue was bounded by $1/b_2\le \lambda_{\min}\le \lambda_{\max}\le 1/b_1$, we have, for sufficiently large $n,$ $\{||\Omega_0-\Omega||_F\le \epsilon_n/b_2\}\subset \mathcal{A}_1.$
Because the
absolute value of the eigenvalues of $\Delta_\Omega$ are bounded by $\epsilon_n/b_2$, for sufficiently large $n$ (i.e. small enough $\epsilon_n$) -- the spectra of $\Omega=\Omega_0+\Delta_\Omega$ is bounded by $\lambda_{\min}-\epsilon_n/b_2$ and $\lambda_{\max}+\epsilon_n/b_2$ and can be further bounded by $\lambda_{\min}/2$ and $2\lambda_{\max}$ when $n$ is large. Thus we can conclude $\tilde{\Pi}(\Omega\succ \tau I|\mathcal{A}_1^\star)=1$ when $n$ is large enough. 

Consequently, we can bound the prior mass $\Pi(\mathcal{A}_1)$ by bounding the prior mass on $\Pi(\mathcal{A}_1^\star)$. Instead of calculating the probability directly, we can lower bound it by observing that
\begin{align*}
    2\sum_{k>k'} |\omega_{0,k,k'}-\omega_{k,k'}|+\sum_k |\omega_{0,k,k}-\omega_{k,k}|2\sum_{(k,k')\in S_0^\Omega} |\omega_{0,k,k'}-\omega_{k,k'}| + 2\sum_{(k,k')\in (S_0^\Omega)^c} |\omega_{k,k'}|+\sum_{k} |\omega_{0,k,k}-\omega_{k,k}|
\end{align*}

Consider the following events:
\begin{align*}
    \mathcal{B}_1=\{ \sum_{(k,k')\in S_0^\Omega} |\omega_{0,k,k'}-\omega_{k,k'}|\le \frac{\epsilon_n}{6b_2}  \},~~\mathcal{B}_2=\{\sum_{(k,k')\in (S_0^\Omega)^c} |\omega_{k,k'}|\le \frac{\epsilon_n}{6b_2}\},~~\mathcal{B}_3&=\{ \sum_{k} |\omega_{0,k,k}-\omega_{k,k}| \le  \frac{\epsilon_n}{3b_2}\} 
\end{align*}
Let $\mathcal{B}=\bigcap_{i=1}^3\mathcal{B}_i\subset \mathcal{A}_1^*\subset \mathcal{A}_1$. Since the prior mass on $\mathcal{B}$ using the lower bound in ~\eqref{eqn:graphical_prior_bound}, we can focus on estimating $\tilde{\Pi}(\mathcal{B})$. Using the fact that $\tilde{\Pi}$ is separable we have
\begin{align*}
    \Pi(\mathcal{A}_1\cap \{\Omega\succ \tau I\})\ge\Pi(\mathcal{A}_1^*\cap \{\Omega\succ \tau I\})\ge \tilde{\Pi}(\mathcal{A}_1^*)\ge \tilde{\Pi}(\mathcal{B})=\prod_{i=1}^3\tilde{\Pi}(\mathcal{B}_i).
\end{align*}

We bound the prior mass on $\mathcal{B}_1$ using pure slab part of the prior, the same technique used in \citet{Bai2020_groupSSL} (specifically Equation D.18). We have
\begin{align*}
\tilde{\Pi}(\mathcal{B}_1)&=\int_{\mathcal{B}_1}\prod_{(k,k')\in S_0^\Omega}   \pi(\omega_{k,k'}|\eta)d\mu\ge \prod_{(k,k')\in S_0^\Omega}  \int_{|\omega_{0,k,k'}-\omega_{k,k'}|\le \frac{\epsilon_n}{6s_0^\Omega b_2}} \pi(\omega_{k,k'}|\eta) d\omega_{k,k'}\\
&\ge \eta^{s_0^\Omega} \prod_{(k,k')\in S_0^\Omega}  \int_{|\omega_{0,k,k'}-\omega_{k,k'}|\le \frac{\epsilon_n}{6s_0^\Omega b_2}} \frac{\xi_1}{2} \exp(-\xi_1|\omega_{k,k'}|) d\omega_{k,k'}\\
&\ge \eta^{s_0^\Omega} \exp(-\xi_1\sum_{(k,k')\in S_0^\Omega}|\omega_{0,k,k'}|)\prod_{(k,k')\in S_0^\Omega}  \int_{|\omega_{0,k,k'}-\omega_{k,k'}|\le \frac{\epsilon_n}{6s_0^\Omega b_2}} \frac{\xi_1}{2} \exp(-\xi_1|\omega_{0,k,k'}-\omega_{k,k'}|) d\omega_{k,k'}\\
&=\eta^{s_0^\Omega}\exp(-\xi_1||\Omega_{0,S_0^\Omega}||_1) \prod_{(k,k')\in S_0^\Omega} \int_{|\Delta| \le \frac{\epsilon_n}{6s^\Omega_0b_2}}\frac{\xi_1}{2} \exp(-\xi_1|\Delta|)d\Delta \\
&\ge \eta^{s_0^\Omega}\exp(-\xi_1||\Omega_{0,S_0^\Omega}||_1)\left[e^{-\frac{\xi_1\epsilon_n}{6b_2s^\Omega_0}}\left(\frac{\xi_1\epsilon_n}{6s^\Omega_0b_2}\right)\right]^{s^\Omega_0}.
\end{align*}

The first inequality holds because the fact that $|\omega_{0,k,k'}-\omega_{k,k'}|\le \epsilon_n/(6s_0^\Omega b_2)$ implies the sum less than $\epsilon_n/(6 b_2)$. The last inequality is a special case of Equation D.18 of \citet{Bai2020_groupSSL}

For $\mathcal{B}_2$ we derive the lower bound using the spike component of the prior. Let $Q=q(q-1)/2$ denote the number of off-diagonal entries of matrix $\Omega$. 
\begin{align*}
\tilde{\Pi}(\mathcal{B}_2)&=\int_{\mathcal{B}_2}\prod_{(k,k')\in (S_0^\Omega)^c}   \pi(\omega_{k,k'}|\eta)d\mu\ge \prod_{(k,k')\in (S_0^\Omega)^c} \int_{|\omega_{k,k'}|\le \frac{\epsilon_n}{6(Q-s_0^\Omega)b_2}} \pi(\omega_{k,k'}|\eta)d\mu\\
&\ge (1-\eta)^{Q-s_0^\Omega} \prod_{(k,k')\in (S_0^\Omega)^c} \int_{|\omega_{k,k'}|\le \frac{\epsilon_n}{6(Q-s_0^\Omega)b_2}} \frac{\xi_0}{2} \exp(-\xi_0|\omega_{k,k'}|)d\omega_{k,k'} \\
& \ge(1-\eta)^{Q-s_0^\Omega}\prod_{(k,k')\in (S_0^\Omega)^c}\left[1-\frac{6(Q-s_0^\Omega)b_2}{\epsilon_n}\mathbb{E}_{\pi}|\omega_{k,k'}|\right]\\
&=(1-\eta)^{Q-s_0^\Omega}\left[1-\frac{6(Q-s_0^\Omega)b_2}{\epsilon_n\xi_0}\right]^{Q-s_0^\Omega}\gtrsim (1-\eta)^{Q-s_0^\Omega} \left[1-\frac{1}{Q-s_0^\Omega}\right]^{Q-s_0^\Omega}\asymp (1-\eta)^{Q-s_0^\Omega}.
\end{align*}

To derive the last line we used the similar argument used by \citet{Bai2020_groupSSL}'s D.22 to establish that $1-\frac{6(Q-s_0^\Omega)b_2}{\epsilon_n\xi_0}\gtrsim 1-1/(Q-s_0^\Omega)$. 
That is, recall our assumption that $\xi_0\sim \max\{Q,n\}^{4+b}$; this assumption allows us to control the $Q$ term and $n$ term induced by $\epsilon_n$ appearing in the numerator in the fifth line.
Because $s_0^\Omega$ grows slower than $Q$, we have $Q-s_0^\Omega\to\infty$ thus we can lower bound the above function using some constant multiple of $(1-\eta)^{Q-s_0^\Omega}$.  

The event $\mathcal{B}_3$ depends only on the diagonal entries of $\Omega$. 
The prior mass can be directly bounded using results on exponential distributions. 
\begin{align*}
\tilde{\Pi}(\mathcal{B}_3)&=\int_{\mathcal{B}_3}\prod_{i=1}^q   \pi(\omega_{k,k})d\mu\ge \prod_{i=1}^q  \int_{|\omega_{0,k,k}-\omega_{k,k}|\le \frac{\epsilon_n}{3q b_2}} \pi(\omega_{k,k}) d\omega_{k,k}=\prod_{i=1}^q  \int_{\omega_{0,k,k}- \frac{\epsilon_n}{3 q b_2}}^{\omega_{0,k,k}+ \frac{\epsilon_n}{3 q b_2}} \xi_1 \exp(-\xi_1\omega_{k,k}) d\omega_{k,k}\\
&\ge \prod_{i=1}^q  \int_{\omega_{0,k,k}}^{\omega_{0,k,k}+ \frac{\epsilon_n}{3 q b_2}} \xi_1 \exp(-\xi_1\omega_{k,k}) d\omega_{k,k}=\exp(-\xi_1 \sum_{i=1}^q\omega_{0,k,k})\int_{0}^{ \frac{\epsilon_n}{3 q b_2}} \xi_1 \exp(-\xi_1\omega_{k,k}) d\omega_{k,k}\\
&\ge \exp(-\xi_1 \sum_{i=1}^q\omega_{0,k,k})\left[e^{-\frac{\xi_1\epsilon_n}{3b_2q}}\left(\frac{\xi_1\epsilon_n}{3qb_2}\right)\right]^{q}
\end{align*}

Now we can show the log prior mass on $\mathcal{B}$ can be lower bounded by some $\exp\{-C_1n\epsilon_n^2\}$, to this end, consider the negative log prior mass
\begin{align}
    \label{eq:B_event_for_Omega}
    \begin{split}
-&\log(\Pi(\mathcal{A}_1\cap \{\Omega\succ \tau I\}))
\le \sum_{i=1}^3 -\log(\tilde{\Pi}(\mathcal{B}_i)) \\
\lesssim& -s_0^\Omega \log(\eta) + \xi_1||\Omega_{0,S_0^\Omega}||_1 + \frac{\xi_1 \epsilon_n}{6b_2}-s_0^\Omega \log\left(\frac{\xi_1\epsilon_n}{6s^\Omega_0b_2}\right)-
(Q-s_0^\Omega) \log(1-\eta)+\xi_1\sum_k \omega_{0,k,k}+\frac{\xi_1\epsilon_n}{3b_2}-q\log\left(\frac{\xi_1\epsilon_n}{3qb_2}  \right)\\
=&-\log\left( \eta^{s_0^\Omega} (1-\eta)^{Q-s_0^\Omega}  \right) + \xi_1||\Omega_{0,S_0^\Omega}||_1 + \frac{\xi_1 \epsilon_n}{6b_2}+\xi_1\sum_k \omega_{0,k,k}+\frac{\xi_1\epsilon_n}{3b_2}-s_0^\Omega \log\left(\frac{\xi_1\epsilon_n}{6s^\Omega_0b_2}\right)-q\log\left(\frac{\xi_1\epsilon_n}{3qb_2}  \right)\\
    \end{split}
\end{align}

The $\frac{\xi_1 \epsilon_n}{6b_2}$ and $\frac{\xi_1\epsilon_n}{3b_2}$ terms are of $O(\epsilon_n)\lesssim n\epsilon_n^2$. The 4th term is of order $q$ since the diagonal entries is controlled by the largest eigenvalue of $\Omega$ that was assumed to be bounded thus $\xi_1\sum_k \omega_{0,k,k}\lesssim n\epsilon_n^2$. Finally, 
\begin{align*}
\xi_1||\Omega_{0,S_0^\Omega}||_1\le \xi_1 s_0^\Omega \sup|\omega_{0,k,k'}|=O(s_0^\Omega)\lesssim n\epsilon_n^2
\end{align*}
because the entries of $\omega_{0,k,k'}$ is bounded. 

Without tuning $\eta$, the first term $-\log\left( \eta^{s_0^\Omega} (1-\eta)^{Q-s_0^\Omega}  \right)$ has order of $Q$.  
But since we assumed $\frac{1-\eta}{\eta}\sim Q^{2+a}$ for some $a>0$, we have $M_1 Q^{2+a} \le \frac{1-\eta}{\eta}\le M_2 Q^{2+a}$. That is, we have $1/(1+M_2 Q^{2+a})\le \eta \le 1/(1+M_1Q^{2+a})$, and a simple lower bound is:
\begin{align*}
\eta^{s_0^\Omega} (1-\eta)^{Q-s_0^\Omega}&\ge (1+M_2 Q^{2+a})^{-s_0^\Omega}(1-\eta)^{Q-s_0^\Omega} \\ 
&\ge (1+M_2 Q^{2+a})^{-s_0^\Omega}\left(1-\frac{1}{1+M_1Q^{2+a}}\right)^{Q-s_0^\Omega}\gtrsim (1+M_2 Q^{2+a})^{-s_0^\Omega}
\end{align*}

The last line is because $Q^{2+a}$ grows faster than $Q-s^\Omega_0$, thus $(1-\frac{1}{1+M_1Q^{2+a}})^{Q-s_0^\Omega}$ can be bounded below with some constant. 
Taking logarithms, we see 
\begin{align*}
-\log\left( \eta^{s_0^\Omega} (1-\eta)^{Q-s_0^\Omega}  \right) \lesssim s_0^\Omega \log(1+M_2 Q^{2+a})\lesssim s_0^\Omega \log(Q)\sim s_0^\Omega\log(q)\le \max(q,s_0^\Omega) \log(q).
\end{align*}

The last two terms in~\eqref{eq:B_event_for_Omega} can be treated in the same way, using our assumption that $\xi_1\asymp1/\max\{n,Q\}$:
\begin{align*}
-s_0^\Omega \log\left(\frac{\xi_1\epsilon_n}{6s^\Omega_0b_2}\right)&=s_0^\Omega \log\left(\frac{6s^\Omega_0b_2}{\xi_1\epsilon_n}\right)\lesssim s_0^\Omega \log\left(\frac{n^{1/2}\max\{n,Q\}s_0^\Omega }{\sqrt{\max\{s_0^\Omega,s_0^B,p,q\}\log(q)}}\right) \\
& \le s_0^\Omega \log\left(n^{1/2}\max\{n,Q\}s_0^\Omega \right)\lesssim s_0^\Omega\log(q^2)\lesssim n\epsilon_n^2
\end{align*}

The third bound holds because $\max\{s_0^\Omega,s_0^B,p,q\}\log(q)\ge 1$ when $n$ is large.
The fourth line follows because(i)  $\log(n)\lesssim \log(q)$ by assumption and (ii) $s_0^\Omega <q^2$. 
The last bound uses the definition of $\epsilon_n$.  

The last term in~\eqref{eq:B_event_for_Omega} can be bounded similar to the one before using our assumption about the slab penalty $\xi_1$:
\begin{align*}
-q \log\left(\frac{\xi_1\epsilon_n}{3qb_2}\right)&=q \log\left(\frac{3qb_2}{\xi_1\epsilon_n}\right)\lesssim q \log\left(\frac{n^{1/2}\max\{n,Q\}q }{\sqrt{\max\{s_0^\Omega,s_0^B,p,q\}\log(q)}}\right) \\ &
\le q \log\left(n^{1/2}\max\{n,Q\}q \right)\lesssim q\log(q)\lesssim n\epsilon_n^2.
\end{align*}
Combining all the results we conclude that $-\log(\Pi(\mathcal{A}_1\cap \{\Omega\succ \tau I\}))\lesssim n\epsilon_n^2$ thus the prior mass $\Pi(\mathcal{A}_1\cap \{\Omega\succ \tau I\})$ is lower bounded as stated in the lemma. 

To bound $\Pi(\mathcal{A}_2|\mathcal{A}_1)$ we used a very similar strategy as the one on $\Omega$. We can also show that mass on L1 norm ball serves as a lower bound similar to that of $\Omega$. Using the definition of KL divergence and variance, the event $\mathcal{A}_2$ can be written as
\begin{align*}
    \mathcal{A}_2=\left\{(B,\Omega): \sum_{i=1}^n||\Omega^{1/2}(B-B_0)^{\top}X_i^{\top}||_2^2\le n\epsilon_n^2, \sum_{i=1}^n||\Omega_0^{-1/2}\Omega(B-B_0)^{\top}X_i^{\top}||_2^2\le n\epsilon_n^2/2\right\}
\end{align*}

Again using the argument from \cite{Ning2020}'s Lemma 5.1. 
Both $n^{-1}\sum_{i=1}^n||\Omega^{1/2}(B-B_0)^{\top}X_i^{\top}||_2^2$ and $n^{-1}\sum_{i=1}^n||\Omega_0^{-1/2}\Omega(B-B_0)^{\top}X_i^{\top}||_2^2$ are bounded by a constant multiplier of 
\begin{align*}
n^{-1}\sum_{i=1}^n||(B-B_0)^{\top}X_i^{\top}||_2^2&=n^{-1}||X(B-B_0)||_F^2.
\end{align*}
We used the facts that $||AB||_F\le \min(||A||_2||B||_F, ||A||_F||B||_2)$ and that $\Omega$ has bounded spectrum.
Since the columns of $X$ have norm of $\sqrt{n}$, we observe that
\begin{align*}
||X\Delta_B||_F\le \sqrt{n}\sum_{j=1}^p||\Delta_{B,j,.}||_F\le \sqrt{n}\sum_{j=1}^p\sum_{k=1}^q|\beta_{j,k}-\beta_{0,j,k}|=\sum_{(j,k)\in S_0^B}|\beta_{j,k}-\beta_{0,j,k}|+\sum_{(j,k)\in (S_0^B)^c}|\beta_{j,k}|
\end{align*}
Thus, to lower bound $\Pi(\mathcal{A}_2|\mathcal{A}_1)$, it suffices to bound $\Pi(\sum |\beta_{j,k}-\beta_{0,j,k}|\le c_1\epsilon_n)$ for some fixed constant $c_1$. We separate the sum based on whether the true value is 0, similar to our treatment of $\Omega$ in the last equality and use the same technique. We consider the event whose intersection is a subset of $\mathcal{A}_2$
\begin{align*}
\mathcal{B}_4=\{ \sum_{(j,k)\in S_0^B}|\beta_{j,k}-\beta_{0,j,k}|\le \frac{c_1\epsilon_n}{2}  \},~~~~
\mathcal{B}_5=\{ \sum_{(j,k)\in (S_0^B)^c}|\beta_{j,k}-\beta_{0,j,k}|\le \frac{c_1\epsilon_n}{2}  \}
\end{align*}
We note $\mathcal{B}_4\cap \mathcal{B}_5\subset\mathcal{A}_2.$
Since the prior for $B$ is separable and independent of $\Omega,$ we have
\begin{align*}
\Pi(\mathcal{A}_2|\mathcal{A}_1)\ge \Pi(\mathcal{B}_4|\mathcal{A}_1)\Pi(\mathcal{B}_5|\mathcal{A}_1)=\Pi(\mathcal{B}_4)\Pi(\mathcal{B}_5).
\end{align*}

We bound each term similarly to $\Omega$. 
To bound $\mathcal{B}_4,$ we use the slab component of $B$'s prior
\begin{align*}
\Pi(\mathcal{B}_4)&=\int_{\mathcal{B}_4}\prod_{(j,k)\in S_0^B}\pi(\beta_{j,k}|\theta)d\mu\ge \prod_{(j,k)\in S_0^B}\int_{|\beta_{j,k}-\beta_{0,j,k}|\le \frac{c_1\epsilon_n}{2s_0^B}}\pi(\beta_{j,k}|\theta)d\beta_{j,k} \\
& \ge \theta^{s_0^B}\prod_{(j,k)\in S_0^B}\int_{|\beta_{j,k}-\beta_{0,j,k}|\le \frac{c_1\epsilon_n}{2s_0^B}}\frac{\lambda_1}{2}\exp(-\lambda_1|\beta_{j,k}|)d\beta_{j,k}\\
&\ge \theta^{s_0^B}\exp(-\lambda_1\sum_{(j,k)\in S_0^B}|\beta_{0,j,k}|)\prod_{(j,k)\in S_0^B}\int_{|\beta_{j,k}-\beta_{0,j,k}|\le \frac{c_1\epsilon_n}{2s_0^B}} \frac{\lambda_1}{2}\exp(-\lambda_1|\beta_{j,k}-\beta_{0,j,k}|)d\beta_{j,k}\\
&=\theta^{s_0^B}\exp(-\lambda_1\sum_{(j,k)\in S_0^B}|\beta_{0,j,k}|)\prod_{(j,k)\in S_0^B}\int_{|\Delta|\le \frac{c_1\epsilon_n}{2s_0^B}} \frac{\lambda_1}{2}\exp(-\lambda_1|\Delta|)d\Delta \\
&\ge \theta^{s_0^B}\exp(-\lambda_1||B_{0,S_0^B}||_1)\left[e^{-\frac{c_1\lambda_1\epsilon_n}{2s_0^B}}\frac{c_1\epsilon_n}{2s_0^B}\right]^{s_0^B}
\end{align*}
In the last line, we used Markov's inequality.
Arguing like before, we can use the spike component of $B$'s prior to control $\P(\mathcal{B}_{5})$:
\begin{align*}
\Pi(\mathcal{B}_5)&\ge (1-\theta)^{pq-s_0^B}\left[1-\frac{2(pq-s_0^B)c_1}{\epsilon_n\lambda_0}\right]^{pq-s_0^B}\gtrsim (1-\theta)^{pq-s_0^B}
\end{align*}
Combining the bounds for $\mathcal{B}_4$ and $\mathcal{B}_5,$ we obtain
\begin{align*}
-\log(\Pi(\mathcal{A}_2|\mathcal{A}_1))&\le -\log(\Pi(\mathcal{B}_4))-\log(\Pi(\mathcal{B}_5))=-\log(\theta^{s_0^B}(1-\theta)^{pq-s_0^B})+\lambda_1||B_{0,S_0^B}||_1+\frac{\lambda_1c_1\epsilon_n}{2}-s_B^0\log\left(\frac{c_1\epsilon_n}{2s_0^B}\right).
\end{align*}

By assumption on the L$_\infty$ norm of $B_0$, we have the entry of $B_0$ is bounded thus the second term $\lambda_1||B_{0,S_0^B}||_1=O(s_0^B)\lesssim n\epsilon_n^2$. The last two terms are $O(\epsilon_n)\lesssim n\epsilon_n^2$.

Recall that we tuned $\frac{1-\theta}{\theta}\sim (pq)^{2+b}$ for some $b>0.$
That is, there are constants $M_{3}$ and $M_{4}$ such that $M_3(pq)^{2+b}\le \frac{1-\theta}{\theta} \le M_4(pq)^{2+b}\le \frac{1-\theta}{\theta}$. 
The resulting bounds on $\theta$ imply that
\begin{align*}
    \theta^{s_0^B}(1-\theta)^{pq-s_0^B}&\ge(1+M_4(pq)^{2+b})^{-s_0^B}(1-\theta)^{pq-s_0^B} \\ 
    &\ge(1+M_4(pq)^{2+b})^{-s_0^B}\left(1-1/(1+M_3(pq)^{2+b})\right)^{pq-s_0^B}\gtrsim (1+M_4(pq)^{2+b})^{-s_0^B}
\end{align*}    
The last bound follows because $(pq)^{2+\beta}$ grows faster than $pq-s_0^B.$ 
Thus $\left(1-1/(1+M_3(pq)^{2+b})\right)^{pq-s_0^B}$ can be bounded from below with some constant not depending on $n$. 
Consequently,
    \begin{align*}
    -\log\left(\theta^{s_0^B}(1-\theta)^{pq-s_0^B}\right)\lesssim s_0^B\log(1+M_4(pq)^{2+b})\lesssim s_0^B\log(pq)\lesssim s_0^B\max\{\log(q),\log(p)\}
    \end{align*}
    To bound the last term we write out $\epsilon_n$ 
\begin{align*}
    -s_B^0\log\left(\frac{c_1\epsilon_n}{2s_0^B}\right)&=s_B^0\log\left(\frac{2s_0^B}{c_1\epsilon_n}\right)\lesssim s_0^B\log\left(\frac{\sqrt{n}s_0^B}{\sqrt{\max\{q,s_0^B,s_0^\Omega\}\log(\max\{p,q\})}}\right)\\
    &\le s_0^B\log(\sqrt{n s_0^B})\lesssim s_0^B\log(pq)\lesssim n\epsilon_n^2
\end{align*}
The third line used the fact when $n$ large $\log(\max\{p,q\})>1$ and the last two lines follow from our assumption that $\log(n)\lesssim \log(q)$ and $s_0^B< pq$. 

Combining all these results, the prior mass $\Pi(\mathcal{A}_2|\mathcal{A}_1)$ can be lower bounded as stated. 
\end{proof}

Lemma~\ref{lemma:dimension_recovery_mssl} guarantees that the posterior places vanishingly small (i.e., $o(1)$) probability on the event that each of $\Omega$ and $B$ have large effective dimensions (i.e., too many entries with large absolute values). Before proceeding to the Lemma, we restate Lemma S3 from \citet{shen2022cgSSL} that allows us to bound probabilities with respect to the prior on $\Omega$ by considering an untruncated version of the prior.
\begin{lemma}
    Let $\tilde{\Pi}$ be the untruncated version of the prior on $\Omega.$
    Then for all events $A,$ for large enough $n$ there is a number $R$ that does not depend on $n$ such that
        \begin{equation}
            \label{eqn:graphical_prior_bound}
            \tilde{\Pi}(\Omega\succ \tau I|A)\tilde{\Pi}(A)\le \Pi_\Omega(A\cap \{\Omega\succ \tau I\})\le \exp(2 \xi_1 Q-\log(R))\tilde{\Pi}(A)
        \end{equation}
        where $Q=q(q-1)/2$ is the total number of free off-diagonal entries in $\Omega.$
\end{lemma}
Now we can state the formal result of dimension recovery of mSSL posterior
\begin{lemma}[Dimension recovery]
\label{lemma:dimension_recovery_mssl}
Letting $s^\star=\max\{q,s_0^\Omega,s_0^B\}, $ for a sufficiently large number $C_3'>0$, we have:
\begin{align}
\sup_{B\in\mathcal{B}_0,\Omega\in\mathcal{H}_0}\E_0\Pi\left(B:\nu_{\beta}(B)>C_3's^\star \vert \bY \right)&\to 0 \label{eq:dimension_B} & &\text{and} &
\sup_{B\in\mathcal{B}_0,\Omega\in\mathcal{H}_0}\E_0\Pi\left(\Omega: \nu_{\omega}(\Omega) >C_3's^\star \vert \bY \right)&\to 0.
\end{align}
\end{lemma}

\begin{proof}
    The main idea is to check the posterior probability directly. 
Let $\mathcal{B}^B_n=\{B:|\nu_{\beta}(B)|<r^B_n\}$ for some $r_n^B=C_3'\max\{q,s_0^B,s_0^\Omega\}$ with $C_3'>C_1$ in the KL condition. 
For $\Omega$, let $\mathcal{B}^\Omega_n=\{\Omega\succ \tau I:|\nu_{\omega}(\Omega)|<r^\Omega_n\}$ for $r_n^\Omega=C_3'\max\{q,s_0^B,s_0^\Omega\}$ with some $C_3'>C_1$ in the KL condition. 
We aim to show that $\E_0\Pi(\Omega\in(\mathcal{B}^\Omega_n)^c|Y_1,\dots,Y_n )\to 0$ and $\E_0\Pi(B\in(\mathcal{B}^B_n)^c|Y_1,\dots,Y_n )\to 0$.

The marginal posterior can be expressed using log-likelihood $\ell_n$:
\begin{equation}
    \label{eqn:posterior_B_m}
    \begin{aligned}
        \Pi(B\in \mathcal{B}_n^B|Y_1,\dots,Y_n)&=\frac{\iint_{\mathcal{B}_n^B}\exp(\ell_n(B,\Omega)-\ell_n(B_0,\Omega_0))d\Pi(B)d\Pi(\Omega)}{\iint\exp(\ell_n(B,\Omega)-\ell_n(B_0,\Omega_0))d\Pi(B)d\Pi(\Omega)}\\
        \Pi(\Omega\in \mathcal{B}_n^\Omega|Y_1,\dots,Y_n)&=\frac{\iint_{\mathcal{B}_n^\Omega}\exp(\ell_n(B,\Omega)-\ell_n(B_0,\Omega_0))d\Pi(B)d\Pi(\Omega)}{\iint\exp(\ell_n(B,\Omega)-\ell_n(B_0,\Omega_0))d\Pi(B)d\Pi(\Omega)}
    \end{aligned}
\end{equation}

By using the result of KL condition (Lemma \ref{lemma:KL_m}), we can, with high probability, bound the denominator from below by $e^{-C_1n\epsilon_n^2}.$
Consequently, focus now on upper bounded the numerators in these expressions, starting with $B.$
\begin{align*}
    \E_0\left(\iint_{(\mathcal{B}_n^{B})^c}f/f_0d\Pi(B)d\Pi(\Omega)\right)
    \le \int_{(\mathcal{B}_n^{B})^c}d\Pi(B)=\Pi(|\nu_{\beta}(B)|\ge r_n^B).
\end{align*}

Notice that when $|\beta_{j,k}|>\delta_\beta,$ we have $\pi(\beta_{j,k})<2\theta\frac{\lambda_1}{2}\exp(-\lambda_1|\beta_{j,k}|).$
Therefore, we can bound the above display as follows
\begin{align*}
    \Pi(|\nu_{\beta}(B)|\ge r_n^B)&\le \sum_{|S|>r_n^B}(2\theta)^{|S|}\prod_{(j,k)\in S}\int_{|\beta_{j,k}|>\delta_\beta} \frac{\lambda_1}{2}\exp(-\lambda_1|\beta_{j,k}|) d\beta_{j,k} \prod_{(j,k)\notin S} \int_{|\beta_{j,k}|<\delta_\beta}\pi(\beta_{j,k})d\beta_{j,k} \\ 
    &\le \sum_{|S|>r_n^B}(2\theta)^{|S|}\\
\end{align*}
Using the assumption on $\theta$, and the fact $\binom{pq}{k}\le (epq/k)^k$ (similar to \citet{Bai2020_groupSSL}'s Equation D.32), we can further upper bound the probability:
\begin{align*}
    \Pi(|\nu_{\beta}(B)|\ge r_n^B)&\le \sum_{|S|>r_n^B}(2\theta)^{|S|}
    \le \sum_{|S|>r_n^B}(\frac{2}{1+M_4 (pq)^{2+b}})^{|S|} \\ 
    &\le \sum_{k=\left\lfloor r_n^B\right\rfloor +1}^{pq} \binom{pq}{k}\left(\frac{2}{M_4(pq)^2}\right)^k
    \le \sum_{k=\left\lfloor r_n^B\right\rfloor +1}^{pq} \left(\frac{2e}{M_4kpq}\right)^k\\
    &<\sum_{k=\left\lfloor r_n^B\right\rfloor +1}^{pq} \left(\frac{2e}{M_4(\left\lfloor r_n^B\right\rfloor +1)pq}\right)^k\lesssim (pq)^{-(\left\lfloor r_n^B\right\rfloor +1)}
    \le \exp(-(\left\lfloor r_n^B\right\rfloor)\log(pq)).
\end{align*}
Taking $r_n^B=C_3'\max\{q,s_0^B,s_0^\Omega\}$ for some $C_3'>C_1$, we have:
\begin{align*}
    \Pi(|\nu_{\beta}(B)|\ge r_n^B)&\le \exp(-C_3'\max\{q,s_0^B,s_0^\Omega\}\log(pq)).
\end{align*}
Therefore
\begin{align*}
    \E_0\Pi((\mathcal{B}_n^B)^c|Y_1,\dots,Y_n)\le \E_0\Pi((\mathcal{B}_n^B)^c|Y_1,\dots,Y_n)I_{E_n}+\P_0(E_n^c),
\end{align*}
where $E_n$ is the event in KL condition. 
On the event $E_n$, the KL condition ensures that the denominator in Equation~\eqref{eqn:posterior_B_m} is lower bounded by $\exp(-C_1n\epsilon_n^2)$ and is upper bounded by $\exp(-C_3'\max\{q,s_0^B,s_0^\Omega\}\log(pq)).$
Further since $\P_0(E_n^c)\to 0$, we conclude that
\begin{align*}
    \E_0\Pi((\mathcal{B}_n^B)^c|Y_1,\dots,Y_n)&\le \exp(C_1n\epsilon_n^2-C_3'\max\{q,s_0^B,s_0^\Omega\}\log(pq))+o(1)\to 0
\end{align*}
The workflow for $\Omega$ is very similar, except we need to use the upper bound of the graphical prior in Equation~\eqref{eqn:graphical_prior_bound} to properly bound the prior mass of the set with overestimated dimensionality of $\Omega$. 

We upper bound the numerator:
\begin{align*}
    \E_0\left(\iint_{(\mathcal{B}_n^{\Omega})^c}f/f_0d\Pi(B)d\Pi(\Omega)\right)&\le \int_{(\mathcal{B}_n^{\Omega})^c}d\Pi(\Omega)=\Pi(|\nu_{\omega}(\Omega)|\ge r_n^\Omega)\le \exp(2\xi_1 Q-\log(R))\tilde{\Pi}(|\nu_{\omega}(\Omega)|\ge r_n^\Omega)
\end{align*}

Whenever $|\omega_{k,k'}|>\delta_\omega,$ we know $\pi(\omega_{k,k'})<2\eta\frac{\xi_1}{2}\exp(-\xi_1|\omega_{k,k'}|)$, yielding the bound 
\begin{align*}
    \tilde{\Pi}(|\nu_{\omega}(\Omega)|\ge r_n^\Omega)&\le \sum_{|S|>r_n^\Omega}(2\eta)^{|S|}\prod_{(k,k')\in S}\int_{|\omega_{k,k'}|>\delta_\omega} \frac{\xi_1}{2}\exp(-\xi_1|\omega_{k,k'}|) d\omega_{k,k'} \prod_{(k,k')\notin S} \int_{|\omega_{k,k'}|<\delta_\omega}\pi(\omega_{k,k'})d\omega_{k,k'} \\
    &\le \sum_{|S|>r_n^\Omega}(2\eta)^{|S|}
\end{align*}

By using the assumption on $\eta$, and the fact $\binom{Q}{k}\le (eQ/k)^k$, we further bound
\begin{align*}
    \tilde{\Pi}(|\nu_{\omega}(\Omega)|\ge r_n^\Omega)&\le \sum_{|S|>r_n^\Omega}(2\eta)^{|S|}
    \le \sum_{|S|>r_n^\Omega}(\frac{2}{1+K_4 Q^{2+b}})^{|S|}\le \sum_{k=\left\lfloor r_n^\Omega\right\rfloor +1}^{Q} \binom{Q}{k}\left(\frac{2}{K_4Q^2}\right)^k \\
    &\le \sum_{k=\left\lfloor r_n^\Omega\right\rfloor +1}^{Q} \left(\frac{2e}{K_4kQ}\right)^k<\sum_{k=\left\lfloor r_n^\Omega\right\rfloor +1}^{Q} \left(\frac{2e}{K_4(\left\lfloor r_n^\Omega\right\rfloor +1)Q}\right)^k\\
    &\lesssim Q^{-(\left\lfloor r_n^\Omega\right\rfloor +1)}
    \le \exp(-(\left\lfloor r_n^\Omega\right\rfloor)\log(Q)).
\end{align*}

Taking $r_n^\Omega=C_3'\max\{q,s_0^B,s_0^\Omega\}$ with $C_3'>C_1$, we have:
\begin{align*}
    \tilde{\Pi}(|\nu_{\omega}(\Omega)|\ge r_n^\Omega)&\le \exp(-C_3'\max\{q,s_0^B,s_0^\Omega\}\log(Q))\le \exp(-C_3'n\epsilon_n^2).
\end{align*}

Using the assumption that $\xi_1 \asymp 1/\max\{Q,n\},$ for some $R'$ not depending on $n$ we conclude that
\begin{align*}
    \Pi(|\nu_{\omega}(\Omega)|\ge r_n^\Omega)\le \exp(-C_3n\epsilon_n^2+2\xi_1 Q-\log(R))\le \exp(-C_3'n\epsilon_n^2+\log(R')).
\end{align*}

Therefore,
\begin{align*}
    \E_0\Pi((\mathcal{B}_n^\Omega)^c|Y_1,\dots,Y_n)\le \E_0\Pi((\mathcal{B}_n^\Omega)^c|Y_1,\dots,Y_n)I_{E_n}+P_0(E_n^c),
\end{align*}
where $E_n$ is the event in KL condition. 
On $E_n$, the KL condition ensures that the denominator in Equation~\eqref{eqn:posterior_B_m} is lower bounded by $\exp(-C_1n\epsilon_n^2)$ and is upper bounded by $\exp(-C_3'n\epsilon_n^2+\log(R')).$ 
Further since $P_0(E_n^c)\to 0$, we obtain the following lower bound
\begin{align}
    \E_0\Pi((\mathcal{B}_n^\Omega)^c|Y_1,\dots,Y_n)&\le \exp(C_1n\epsilon_n^2-C_3'n\epsilon_n^2+\log(R'))+o(1)\to 0
    \label{eq:Omega_dimension}
\end{align}
\end{proof}

Lemma~\ref{lemma:dimension_recovery_mssl} allows us to focus our attention on those $B$'s and $\Omega$'s that do not have too many non-zero entries.
Specifically, we define
\begin{align*}
\mathcal{B}^B_n&:=\left\{B:|\nu_{\beta}(B)|<C_3'\max\{q,s_0^B,s_0^\Omega\}\right\} \\
\mathcal{B}^\Omega_n &:=\left\{\Omega\succ \tau I:|\nu_{\omega}(\Omega)|<C_3'\max\{q,s_0^B,s_0^\Omega\}\right\}.
\end{align*}
Then, we defined the sieve
\begin{equation}
\label{eqn:sieve}
\mathcal{F}_n:=\left\{B\in \mathcal{B}^B_n,\Omega\in \mathcal{B}^\Omega_n:||B||_1\le 2C_3p, ||\Omega||_1 \le 8C_3q \right\},
\end{equation}
for some constant $C'_{3} > C_1$ and $C_3 > C_1+2+\log(3)$ not depending on $n$. 
Lemma~\ref{lemma:sieve} shows that the modified mSSL prior places vanishingly small prior probability outside the sieve. 
\begin{lemma}[Sieve]
    \label{lemma:sieve}
$\mathcal{F}_n$ receive overwhelming prior probability: $\Pi(\mathcal{F}_n^c)\le e^{-C_2n\epsilon_n^2}$.
\end{lemma}
\begin{proof}
    We show the first part of Lemma 3 that 
\begin{equation}
    \mathcal{F}_n=\left\{B\in \mathcal{B}^B_n,\Omega\in \mathcal{B}^\Omega_n:||B||_1\le 2C_3p, ||\Omega||_1 \le 8C_3q \right\}
\end{equation}
for some large $C_3>C_1+2+\log(3)$ where $C_1$ is the constant in KL condition. 
We have
\begin{align*}
    \Pi(\mathcal{F}_n^c)&\le \Pi(||B||_1> 2C_3p)+\Pi((||\Omega||_1 > 8C_3q)\cap \{\Omega\succ \tau I\}).
\end{align*}

We upper bound each term. 
Using the bound in Equation~\eqref{eqn:graphical_prior_bound}, we have
\begin{align*}
    \Pi((||\Omega||_1 > 8C_3q)\cap \{\Omega\succ \tau I\})\le \exp(2\xi_1 Q-\log(R))\tilde{\Pi}(||\Omega||_1 > 8C_3q).
\end{align*}

Since $||\Omega||_1=2\sum_{k>k'}|\omega_{k,k'}|+\sum_{k}|\omega_{k,k}|$, at least one of these two sums exceeds $8C_3q/2$. Thus, we can form an upper bound on the L1 norm probability
\begin{align*}
    \tilde{\Pi}\left(||\Omega||_1 > 8C_3q\right)\le \tilde{\Pi}\left(\sum_{k>k'}|\omega_{k,k'}|>\frac{8C_3q}{4}\right)+\tilde{\Pi}\left(\sum_{k}|\omega_{k,k}|>\frac{8C_3q}{2}\right).
\end{align*}

Under $\tilde{\Pi}$, the slab distribution provides an an upper bound, where $\sum_{k>k'}|\omega_{k,k'}|$ is Gamma distributed with shape parameter $Q$ and rate parameter $\xi_1$. 
By using an appropriate tail estimation of Gamma distribution (\citet{boucheron2013concentration}, pp.29), and the fact $1+x-\sqrt{1+2x}\ge (x-1)/2,$ we compute
\begin{align*}
    \exp(2\xi_1 Q-\log(R))\tilde{\Pi}(\sum_{k>k'}|\omega_{k,k'}| > 8C_3q/4)&\le \exp\left[-Q\left(1-\sqrt{1+2\frac{8C_3q}{4Q\xi_1}}+\frac{8C_3q}{4Q\xi_1}\right)+2Q-\log(R)\right]\\
    &\le \exp\left[-\frac{8C_3q}{8\xi_1}+\left(\frac{5}{2}Q-\log(R)\right)\right].
\end{align*}

Using the fact $\xi_1\asymp 1/\max\{n,Q\}$, we have as $n$ sufficiently large, $ n\epsilon_n^2\ge q\log(q)$ thus $qn\epsilon_n^2\ge Q\log(q)$ and $Q=o(qn\epsilon_n^2)$, 
\begin{align*}
    \frac{8C_3q}{8\xi_1}-\left(\frac{5}{2}Q-\log(R)\right)&\asymp C_3(\max\{n,Q\}q)-\left(\frac{5}{2}Q-\log(R)\right)\\
    &\ge C_3(qn\epsilon_n^2)-\left(\frac{5}{2}Q-\log(R)\right)=C_3(qn\epsilon_n^2)-o(qn\epsilon_n^2)\ge C_3n\epsilon_n^2.
\end{align*}
The first order term of $Q$ on the left hand side can be ignored when $n$ large as the left hand side is dominated by the $Q\log(q)$ term. We have also used the assumption that $\epsilon_n\to 0$. 
Thus we have,
\begin{align*}
    \exp(2 \xi_1 Q-\log(R))\tilde{\Pi}((\sum_{k>k'}|\omega_{k,k'}| > 8C_3q/4))\le \exp(-C_3n\epsilon_n^2).
\end{align*}

For the diagonal, the sum is gamma distributed with shape $q$ and rate $\xi_1$, yielding a similar bound:
\begin{align*}
    \exp(2 \xi_1 Q-\log(R))\tilde{\Pi}(\sum_{k}|\omega_{k,k}| > 8C_3q/2)&\le \exp(2 Q-\log(R))\exp\left[-q\left(1-\sqrt{1+2\frac{8C_3q}{2q\xi_1}}+\frac{8C_3q}{2q\xi_1}\right)\right]\\
    &\le \exp\left[-\frac{8C_3q}{4\xi_1}+Q\left(2+\frac{q}{2Q}\right)-\log(R)\right].
\end{align*}

By using the same argument as before and the fact that $\xi_1 \asymp 1/\max\{Q,n\},$ we have
\begin{align*}
\frac{8C_3q}{4\xi_1}-Q\left(2+\frac{q}{2Q}\right)+\log(R)&\asymp 2C_3(\max\{Q,n\}q)-Q\left(2+\frac{q}{2Q}\right)+\log(R)\ge  C_3qn\epsilon_n^2-o(qn\epsilon_n^2)\ge C_3n\epsilon_n^2.
\end{align*}

The first order term of $Q$ on the left hand side can be ignored when $n$ large as the left hand side is dominated by the $Q\log(q)$ term and $q/Q\to 0$. 
By combining the above results, we have 
\begin{equation}
\label{eqn:sieve_Omega_bound_mssl}
\begin{aligned}
	&\Pi((||\Omega||_1 > 8C_3q)\cap \{\Omega\succ \tau I\})\le \exp(2 Q-\log(R))\tilde{\Pi}(||\Omega||_1 > 8C_3q)\\
	\le& \exp(2 Q-\log(R))\tilde{\Pi}(\sum_{k>k'}|\omega_{k,k'}|>\frac{8C_3q}{4})+\exp(2 Q-\log(R))\tilde{\Pi}(\sum_{k}|\omega_{k,k}|>\frac{8C_3q}{2})
	\le 2\exp(-C_3n\epsilon_n^2).
\end{aligned}
\end{equation}
The probability $||B||_1>2C_3p$ can be bound by tail probability of Gamma distribution with shape parameter $pq$ and rate parameter $\lambda_1$:
\begin{align*}
	\Pi(||B||_1> 2C_3p)&\le \exp\left[-pq\left(1-\sqrt{1+2\frac{2C_3p}{pq\lambda_1}}+\frac{2C_3p}{pq\lambda_1}\right)\right]\le  \exp\left[-pq\left(\frac{2C_3p}{2pq\lambda_1}-\frac{1}{2}\right)\right]\le \exp\left(-\frac{2C_3p}{2\lambda_1}+\frac{pq}{2}\right).
\end{align*}
Using the same argument, we have, $pn\ge pn\epsilon_n^2\ge pq\log(q)$ so $pq=o(pn\epsilon_n^2)$ for large $n$, and
\begin{align*}
	\exp\left(-\frac{2C_3p}{2\lambda_1}+\frac{pq}{2}\right)\le \exp\left(-C_3pn\epsilon_n^2+o(pn\epsilon_n^2)\right)\le \exp(-C_3n\epsilon_n^2).
\end{align*}
Combining these two inequalities, we have
\begin{equation}
\label{eqn:sieve_B_bound_mssl}
\begin{aligned}
	\Pi(||B||_1> 2C_3p)&\le\exp\left(-\frac{2C_3p}{2\lambda_1}+\frac{pq}{2}\right) \le\exp(-C_3n\epsilon_n^2)
\end{aligned}
\end{equation}

By combining the result from Equations~\eqref{eqn:sieve_Omega_bound_mssl} and~\eqref{eqn:sieve_B_bound_mssl}, we have
\begin{align*}
	\Pi(\mathcal{F}_n^c)\le 3\exp(-C_3n\epsilon_n^2)=\exp(-C_3n\epsilon_n^2+\log(3)).
\end{align*}
With our choice of $C_3$, the above probability is asymptotically bounded from above by $\exp(-C_2n\epsilon_n^2)$ for some $C_2\ge C_1+2$. 
\end{proof}

Lemma~\ref{lemma:test-condition} shows that on $\mathcal{F}_{n},$ we can construct hypothesis tests with error rates that vanish fast sufficiently fast.
To do so, we first construct Neyman-Pearson tests with vanishing error rates within small-norm balls in the parameter space.
Then, by showing that the number of such balls needed to cover $\mathcal{F}_{n}$ doe not grow too quickly with $n,$ we can construct the required test by taking a supremum over these balls.
A key step in our proof of Lemma~\ref{lemma:test-condition} was bounding the packing number of $\mathcal{F}_{n}$ by packing ``effectively'' low-dimensional sets.
We restate this result, which originally appeared as Lemma S4 in \cite{shen2022cgSSL}, here as Lemma~\ref{lemma:packing_number_lp}.

\begin{lemma}[Lemma S4 in \cite{shen2022cgSSL}]
    \label{lemma:packing_number_lp}
    For a set of form $E=A \times [-\delta,\delta]^{Q-s}\subset \mathbb{R}^Q$ where $A\subset \mathbb{R}^s$, (with $s>0$ and $Q\ge s+1$ are integers) for $1\le p<\infty$ and a given $T>1$, if $\delta<\frac{\epsilon}{2[T(Q-s)]^{1/p}}$, we have the packing number:
    \begin{align*}
        D(\epsilon, A, ||\cdot||_p) \le D(\epsilon,E, ||\cdot ||_p)\le D((1-T^{-1})^{1/p}\epsilon, A, ||\cdot||_p).
    \end{align*}
\end{lemma}

\begin{lemma}[Test condition]
\label{lemma:test-condition}
There exists tests $\varphi_{n}$ with Type I and Type II error bounded by $e^{-M_{2}n\epsilon_{n}^{2}}$ for some constant $M_{2} > C_{1}+1,$ where $C_{1}$ is the constant from Lemma~\ref{lemma:KL_m}.
That is, $\E_{f_{0}}\varphi_{n} \lesssim e^{-M_{2}n\epsilon_{n}^{2}/2}$ and
$$
\sup_{f \in \mathcal{F}_{n}:\rho(f_{0},f) > M_{2}n\epsilon_{n}^{2}}\E_{f}(1-\varphi_{n}) \lesssim e^{-M_2n\epsilon_n^2}.
$$
\end{lemma}

\begin{proof}
    Instead of directly constructing the $\varphi_n$ whose alternative is the whole sieve, we start by constructing a sequence of tests whose alternative correspond tp representative points. 
We show that these tests have bounded Type II error in the neighborhood of the representative points.
We can then construct a test against the entire sieve by taking the supremum of the tests constructed against representative points. 
We finally apply the general theory by bounding the number of neighborhoods around representative points needed to cover the entire sieve. 

For a representative point $f_1$, we consider the Neyman-Pearson test against a single point alternative $H_0: f=f_0, H_1: f=f_1$, $\phi_n=I\{ f_1/f_0\ge 1 \}$. 
If the average half-order R\'enyi divergence $-n^{-1}\log(\int\sqrt{f_0f_1}d\mu)\ge \epsilon^2$,  we will have:
\begin{align*}
    \mathbb{E}_{f_0}(\phi_n)&\le\int_{f_1>f_0}\sqrt{f_1/f_0} f_0 d\mu\le \int \sqrt{f_1f_0}d\mu\le e^{-n\epsilon^2} \\
    \mathbb{E}_{f_1}(1-\phi_n) &\le\int_{f_0>f_1}\sqrt{f_0/f_1} f_1 d\mu\le \int \sqrt{f_0f_1}d\mu\le e^{-n\epsilon^2}.
\end{align*}
By Cauchy-Schwarz, for any alternative $f$ we can control the type-II error rate:
\begin{align*}
    \E_f(1-\phi_n)\le \{\E_{f_1}(1-\phi_n)\}^{1/2}\{\E_{f_1}(f/f_1)^2\}^{1/2}.
\end{align*}
So long as the second factor grows at most like $e^{cn\epsilon^2}$ for some properly chosen small $c$, the full expression can be controlled.  
Thus we can consider the neighborhood around the representative point small enough so that the second factor can be actually bounded.

Consider every density with parameters satisfying
\begin{align}
    \label{eqn:test_sets_mssl}
    \begin{split}
    |||\Omega|||_2\le ||\Omega||_1 & \le 8C_3q, \\
    ||B_1-B||_2\le ||B_1-B||_1 &\le \frac{1}{\sqrt{2C_3npq}} \\
    |||\Omega_1-\Omega|||_2\le ||\Omega_1-\Omega||_1 &\le \frac{1}{8C_3nq^{3/2}}.
\end{split}
\end{align}

We show that $\E_{f_1}(f/f_1)^2$ is bounded on the above set when parameters are from the sieve $\mathcal{F}_n$. 
Denote $\Sigma_1=\Omega_1^{-1}$, $\Sigma=\Omega^{-1}$ as well as $\Sigma_1^\star=\Omega^{1/2}\Sigma_1\Omega^{1/2}$, and $\Delta_B=B-B_1$ while $\Delta_\Omega=\Omega-\Omega_1$. 
\begin{align}
    \begin{split}
    \label{eqn:test_f_over_f1_mssl}
        \E_{f_1}(f/f_1)^2=&|\Sigma_1^\star|^{n/2}|2I-\Sigma_1^{\star-1}|^{-n/2} \times 
        \\ &~~\exp\left(\sum_{i=1}^nX_i(B-B_1)\Omega^{1/2}(2\Sigma_1^\star-I)^{-1}\Omega^{1/2}(B-B_1)^{\top}X_i^{\top}\right).
    \end{split}
\end{align}
For the first factor, 
since $\Omega\in\mathcal{F}_n$, we have $|||\Omega^{-1}|||_2\le 1/\tau.$
The fact that $|||\Omega_1-\Omega|||_2\le \delta_n'= 1/8C_3nq^{3/2}$ implies that
\begin{align*}
    |||\Sigma_1^\star-I|||_2\le |||\Omega^{-1}|||_2|||\Omega_1-\Omega|||_2\le \delta_n'/\tau.
\end{align*}
We can therefore bound the spectrum of $\Sigma_1^\star$, i.e. $1-\delta_n'/\tau\le \eig_1(\Sigma_1^\star)\le \eig_q(\Sigma_1^\star)\le  1+\delta_n'/\tau $.

Thus
\begin{align*}
    \left(\frac{|\Sigma_1^\star|}{|2I-\Sigma_1^{\star-1}|}\right)^{n/2}&=\exp\left(\frac{n}{2}\sum_{i=1}^q\log(\eig_i(\Sigma_1^\star))-\frac{n}{2}\sum_{i=1}^q \log\left(2-\frac{1}{\eig_i(\Sigma_1^\star)}\right)\right) \\
    & \le \exp\left(\frac{nq}{2}\log(1+\delta_n'/\tau)-\frac{nq}{2}\log\left(1-\frac{\delta_n'/\tau}{1-\delta_n'/\tau}\right)\right)\\
    &\le \exp\left(\frac{nq^2}{2}\delta_n'+\frac{nq}{2}\left(\frac{\delta_n'/\tau}{1-2\delta_n'/\tau}\right)\right)\le \exp(nq\delta_n'/\tau)\le e
\end{align*}
The third inequality follows from the fact $1-x^{-1}\le \log(x)\le x-1$.

We can bound the log of the second factor by
\begin{align*}
    |||\Omega|||_2|||(2\Sigma_1^\star-I)^{-1}|||_2\sum_{i=1}^n||X_i\Delta_B||_2^2\le 16C_3q \sum_{i=1}^n||X_i\Delta_B||_2^2\le 16C_3npq ||\Delta_B||_2^2\le 8
\end{align*}
The desired test $\varphi_n$ in Lemma \ref{lemma:test-condition} can be obtained as the maximum of all tests $\phi_n$ described above -- that is, the maximum over all sets defined in Equation~\eqref{eqn:test_sets_mssl} that covers the sieve $\mathcal{F}_n$.

To finish the proof, we use a covering arguments to show that $\phi_n$'s we take maximum over is not too big. Formally we show that number of sets described in Equation~\eqref{eqn:test_sets_mssl} needed to cover sieve $\mathcal{F}_n$, denoted by $N_*$, can be bounded by $\exp(Cn\epsilon_n^2)$ for some suitable constant $C$. 

We can bound the logarithm of the covering number $\log(N_*)$:
\begin{align*}
    \log(N_*)\le & \log\left[N\left(\frac{1}{\sqrt{2C_3npq}},\{B\in \mathcal{B}^B_n:||B||_1\le2C_3p\},||\cdot||_1\right)\right] \\
    & ~~ +\log\left[N\left(\frac{1}{8C_3nq^{3/2}},\{\Omega \in \mathcal{B}^\Omega_n,||\Omega||_1\le 8C_3q\},||\cdot||_1\right)\right]
\end{align*}
The two terms above can be treated in a similar fashions. 
Denote $\max\{q,s_0^B,s_0^\Omega\}=s^\star$. 
Since there are multiple ways to allocate the effective 0's, we pick up a binomial coefficient below:
\begin{equation}   
\begin{aligned}
    &N\left(\frac{1}{8C_3n q^{3/2}},\{\Omega \in \mathcal{B}^\Omega_n,||\Omega||_1\le 8C_3q\},||\cdot||_1\right)\\
    \le & \binom{Q}{C_3's^\star} N\left(\frac{1}{8C_3nq^{3/2}},\{V\in \mathbb{R}^{Q+q}:|v_i|<\delta_\omega\text { for $1\le i\le Q+q-C_3's^\star$},||V||_1\le 8C_3q\},||\cdot||_1\right) \label{eq:covering_Omega} \\
    &N\left(\frac{1}{\sqrt{2C_3npq}},\{B\in \mathcal{B}^B_n:||B||_1\le2C_3p\},||\cdot||_1\right)\\
    \le & \binom{pq}{C_3' s^\star}  N\left(\frac{1}{\sqrt{2C_3npq}},\{V\in \mathbb{R}^{pq}:|v_i|<\delta_\beta\text { for $1\le i\le pq-C_3's^\star$},||V||_1\le 2C_3p\},||\cdot||_1\right).
\end{aligned}
\end{equation}

Note that $\Omega$ has $Q+q<2Q$  free parameters. 
We have first
\begin{align*}
    \log \binom{Q}{C_3's^\star}&\lesssim s^\star \log(Q)\lesssim n\epsilon_n^2,~~~~\log \binom{pq}{C_3's^\star}\lesssim s^\star \log(pq)\lesssim n\epsilon_n^2. 
\end{align*}

Observe that for $V' \in \R^{C_{3}'s^{\star}}$, we have
$$||V||_1\cap \{|v_i|<\delta_\omega \text{ for } 1\le i\le Q+q-C_3's^\star\}\subset \{||V'||\le 8C_3q\}\times[-\delta_\omega,\delta_\Omega]^{Q+q-C_3's^\star}.
$$
Thus,
\begin{align*}
    &N\left(\frac{1}{8C_3nq^{3/2}},\{V:|v_i|<\delta_\omega\text { for $1\le i\le Q+q-C_3's^\star$},||V||_1\le 8C_3q\},||\cdot||_1\right)\\
    \le & N\left(\frac{1}{8C_3nq^{3/2}},\{V\in \mathbb{R}^{C_3's^\star}:||V'||_1\le 8C_3q\times [-\delta_\omega,\delta_\omega]^{Q+q-C_3's^\star}\},||\cdot||_1\right)\\
\end{align*}

Because we are covering the sieve with $L_{1}$-norm balls, we will verify the condition of Lemma 5 with $p = 1$ and $T = 2.$
The choice $T=2$ is arbitrary but makes the calculation easy. 
By our assumption on $\xi_0$, we have
\begin{align*}
    (Q+q-C_3's^\star)\delta_\omega&\le 2Q\delta_\omega = 2Q\frac{1}{\xi_0-\xi_1}\log\left[\frac{1-\eta}{\eta}\frac{\xi_0}{\xi_1}\right]\lesssim  \frac{Q\log(\max\{q,n\})}{\max\{Q,n\}^{4+b/2+b/2}}\le \frac{1}{\max\{Q,n\}^{3+b/2}}
\end{align*}
The denominator dominates $ C_3nq^{3/2}$ thus for large enough $n$, we have $(Q+q-C_3's^\star)\delta_\omega\le \frac{1}{32C_3nq^{3/2}}$ thus by Lemma 5, we can control the covering number by the packing number:
\begin{align*}
    &
    \log N\left(\frac{1}{8C_3nq^{3/2}},\{V:|v_i|<\delta_\omega\text { for $1\le i\le Q+q-C_3's^\star$},||V||_1\le 8C_3q\},||\cdot||_1\right)\\
    \le &\log D\left(\frac{1}{16C_3nq^{3/2}},\{V'\in \mathbb{R}^{C_3's^\star},||V'||_1\le 8C_3q\},||\cdot||_1\right)\\
    \lesssim & s^\star \log(128C_3^2qnq^{3/2})
    \lesssim n\epsilon_n^2.
\end{align*}

We treat the covering number for $B$ in a similar fashion:
\begin{align*}
    &N\left(\frac{1}{\sqrt{2C_3npq}},\{V:|v_i|<\delta_\beta\text { for $1\le i\le pq-C_3's^\star$},||V||_1\le 2C_3p\},||\cdot||_1\right)\\
    \le& N\left(\frac{1}{\sqrt{2C_3npq}},\{V'\in \mathbb{R}^{C_3's^\star}:||V'||_1\le 2C_3p\times [-\delta_\beta,\delta_\beta]^{pq-C_3's^\star}\},||\cdot||_1\right)
\end{align*}

We once again verify the condition of Lemma 5 with $p=1$ and $T=2$:
\begin{align*}
    (pq-C_3's^\star )\delta_\beta&\le pq\delta_\beta=\frac{pq}{\lambda_0-\lambda_1}\log\left[\frac{1-\theta}{\theta}\frac{\lambda_0}{\lambda_1}\right]\lesssim \frac{pq\log(\max\{p,q,n\})}{\max\{pq,n\}^{5/2+b/2+b/2}}\le \frac{1}{\max\{pq,n\}^{3/2+b/2}}.
\end{align*}
The denominator dominates $\sqrt{2C_3npq}$, Thus for enough large $n$, we have $(pq-C_3's^\star )\delta_\beta\le 1/4\sqrt{2C_3npq}$. Thus similar to $\Omega$, we have:
\begin{align*}
    &
    \log N\left(\frac{1}{\sqrt{2C_3npq}},\{V:|v_i|<\delta_\omega\text { for $1\le i\le pq-C_3's^\star$},||V||_1\le 2C_3p\},||\cdot||_1\right)\\
    \le &\log D\left(\frac{1}{2\sqrt{2C_3npq}},\{V'\in \mathbb{R}^{C_3's^\star},||V'||_1\le 2C_3p\},||\cdot||_1\right)
    \lesssim  s^\star \log(4C_3p\sqrt{2C_3npq})
    \lesssim n\epsilon_n^2.
\end{align*}
\end{proof}

As a direct consequence of contraction in log-affinity, we can show that the posterior distribution of $\Omega$ and $XB,$ respectively, concentrate around $\Omega_{0}$ and $XB_{0}.$ 
\begin{theorem}[Posterior contraction of mSSL]
    \label{thm:posterior_contraction_mssl}
    Under Assumptions A1--A5, there is some constant $M_1>0$ that does not depend on $n$ such that
    \begin{align}
        \sup_{B\in\mathcal{B}_0,\Omega\in\mathcal{H}_0}\E_0 \Pi\left(B:||X(B-B_0)||_F^2\ge M_1n\epsilon_n^2|Y_1,\dots,Y_n\right) &\longrightarrow 0\label{eqn:contract_predict_m}\\
        \sup_{B\in\mathcal{B}_0,\Omega\in\mathcal{H}_0}\E_0 \Pi\left(\Omega:||\Omega-\Omega_0||_F^2\ge M_1\epsilon_n^2|Y_1,\dots,Y_n\right) &\longrightarrow 0\label{eqn:contract_omega_m}
    \end{align}
    where $\epsilon_n=\sqrt{\max\{q,s_0^\Omega,s_0^B\}\log(\max\{p,q\})/n}.$
\end{theorem}

\begin{proof}
Since $\sum \rho(f_i-f_{0i})\lesssim n\epsilon_n^2,$ we know
\begin{equation}
    \label{eqn:log_affinity_implies_some_bound_mssl}
    \begin{aligned}
        -\log\left(\frac{\lvert \Omega^{-1} \rvert^{1/4}\lvert\Omega_0^{-1}\rvert^{1/4}}{\vert(\Omega^{-1}+\Omega_0^{-1})/2\rvert^{1/2}}\right)&\lesssim \epsilon_n^2\\
        \frac{1}{8n}\sum \bx_i^{\top}(B - B_{0})\left(\frac{\Omega^{-1}+\Omega_0^{-1}}{2}\right)^{-1}(B - B_{0})^{\top}\bx_i&\lesssim \epsilon_n^2\\
    \end{aligned}
\end{equation}
Because $\Omega^{-1}$ has bounded operator norm (Assumption A1), the first line in Equation~\ref{eqn:log_affinity_implies_some_bound_mssl} implies that $\lVert \Omega^{-1}-\Omega_0^{-1}\rVert_{F}^{2}\lesssim \epsilon_n^2$ \cite[c.f.][Equation 5.11]{Ning2020}.

Meanwhile since $\Omega_{0}$ has bounded spectrum (Assumption A1), the fact that $\lVert \Omega^{-1}-\Omega_0^{-1}\rVert_{F}^{2}\lesssim \epsilon_n^2$ implies that, for large enough $n,$ $\Omega$'s L2 operator norm is bounded.
Using the fact that (i) $\lVert AB\rVert_F\le |||A|||_2 \lVert B\rVert_F$ for matrices $A$ and $B;$ (ii) the decomposition $\Omega_0-\Omega=\Omega(\Omega^{-1}-\Omega_0^{-1})\Omega_0;$ and (iii) Assumption A1 again, we conclude that~\eqref{eqn:log_affinity_implies_some_bound_mssl} implies $\lVert \Omega-\Omega_0\rVert_F\lesssim \epsilon_n.$ 

To prove \eqref{eqn:contract_predict_m}, note that because $|||\Omega^{-1}|||_{2}$ is bounded for large enough $n,$ the second inequality in~\eqref{eqn:log_affinity_implies_some_bound_mssl} implies that
\begin{align*}
    \epsilon_n^2
    &\gtrsim \frac{1}{8n}\sum \frac{||\bx_{i}^{\top}(B - B_{0})||_2^2}{|||(\Omega^{-1} + \Omega_{0}^{-1})/2|||_2}\\
    &\gtrsim\frac{1}{n}\sum \lVert\bx_{i}^{\top}(B - B_{0})\rVert_2^2/\sqrt{\epsilon_n^2+1}.
\end{align*}
Note that \citet{Ning2020} utilizes a similar argument (cf. their Equation 5.12). 
\end{proof}

The result in Equation~\eqref{eqn:contract_predict_m} shows that the matrix $XB$ contracts to its true value $XB_0$ in Frobenius norm. 
Equation~\eqref{eqn:contract_omega_m} guarantees that the residual precision matrix $\Omega$ also contracts to its true value $\Omega_0$. 
Importantly, apart from Assumption A2 about the dimension of $X$, Theorem~\ref{thm:posterior_contraction_mssl} does not require any additional assumptions about the design matrix.
The contraction rate for the matrix of marginal effects $B,$ on the other hand, depends critically on $X,$ through a particular \textit{restricted eigenvalue} defined to be
$$
\phi^{2}(s) = \inf_{\substack{A\in \mathbb{R}^{p\times q}: \\ 0\le |\nu(A)|\le s}} \left\{ \frac{\lVert XA \rVert^{2}_{F}}{n\lVert A \rVert^{2}_{F}} \right\}.
$$

\begin{corollary}[Recovery of regression coefficients in mSSL]
\label{cor:B_Psi_recovery}
    Under Assumptions A1--A5, there is some constant $M' > 0,$ which does not depend on $n,$ such that
     \begin{align}
        \sup_{B\in \mathcal{B}_0,\Omega\in\mathcal{H}_0}\E_0 \Pi\left(||B-B_0||_F^2\ge \frac{M'\epsilon_n^2}{\phi^2(C_3's^\star)}\right) &\to 0 \label{eqn:recover_B_m} 
    \end{align}
where $s^\star=\max\{q,s_0^\Omega,s_0^B\}$.
\end{corollary}

\begin{proof}

We have by dimension result~\eqref{eq:dimension_B} $B$ asymptotically have dimension less than $C_{3}'s^{\star}$, thus 
$$
\lVert XB - XB_{0} \rVert^{2}_{F} \leq \phi^{2}(C_{3}'s^{\star})\lVert B - B_{0} \rVert_{F}^{2}.
$$
Thus, the result in~\eqref{eqn:contract_predict_m} implies that
$$
 \sup_{B\in \mathcal{B}_0,\Omega\in\mathcal{H}_0}\E_0 \Pi\left(||B-B_0||_F^2\ge \frac{M'\epsilon_n^2}{\phi^2(C_3's^\star)}\right)\to 0.
$$

\end{proof}

\subsection{Comparison with known results}
\label{sec:comparisons}

Theorem~\ref{thm:posterior_contraction_mssl} and Corollary~\ref{cor:B_Psi_recovery} generalize several existing posterior contraction results about spike-and-slab posteriors in simpler models.
We briefly summarize these results and describe how they arise as special cases of our results.

\textbf{Single-outcome regression}. In the single-outcome (i.e., $q = 1$) setting with known residual variance and at least one non-zero parameter, \citet{RockovaGeorge2018_ssl} established a posterior contraction rate of $\sqrt{s_{0}^{B}\log{p}/n}$ for  average regression function $XB/n$ and $\sqrt{s_{0}^{B}\log{p}/n\phi^2}$ for regression coefficients under their definition of restricted eigenvalue of $\phi$.
We exactly obtain their rate by plugging $q = 1$ into our formula for $\epsilon_{n}.$

\textbf{Gaussian graphical model}. The Gaussian graphical model is a special case of the model in Equation~\eqref{eq:general_model} with $p = 0$ covariates that assumes $\by \vert \Omega \sim \mathcal{N}_{q}(0,\Omega^{-1}).$
\citet{Gan2019_unequal} studied the sparse Gaussian graphical model with spike-and-slab priors on the off-diagonal elements $\omega_{k,k'}.$
Their MAP estimator had a convergence rate equivalent to $\sqrt{\max\{s_0^\Omega,q\}\log(q)/n}$ in Frobenius norm, which agrees numerically with our $\epsilon_{n}$ when we set $p = 0.$
In fact, with slight modifications, our proof of Theorem~\ref{thm:posterior_contraction_mssl} can be used to establish the posterior contraction rate for the sparse Gaussian graphical model with SSL priors; see Section S2 in the Supplemental Materials.
We note that \citet{banerjee2021precision} obtained the same rates by using horseshoe-like shrinkage priors instead of two-group mixture priors. 

\textbf{Separate regressions.} Recall that $k^{\text{th}}$ column of $B$ contains the marginal effect of each predictor on the outcome $Y_{k}.$
The mSSL leverages information about all of the other outcomes in order to estimate this column.
An alternative strategy, which does not borrow strength across outcomes, is to estimate the columns of $B$ one at a time using \citet{RockovaGeorge2018_ssl}'s SSL for each individual response.
\citet{Deshpande2019} showed empirically that the alternative ``separate SSL'' strategy was not as capable as the mSSL at recovering the support of $B$ or estimating $B.$

When $\Omega$ is truly diagonal, the separate SSL posterior contracts at the rate $\sqrt{s_0^B\log(\max\{p,q\})/n}$ for estimating $B$ in Frobenius norm. 
The rate implied by our results on the mSSL in this special setting is $\sqrt{\max\{s_0^B,q\}\log(\max\{p,q\})/n}.$
So long as there is at least one covariate $X_{j}$ that is predictive of each outcome $Y_{k},$ we know that $s_{0}^{B} \geq q.$ 
In other words, although the mSSL may have better finite-sample empirical performance than the separate SSL, the two approaches are have the same asymptotic convergence rates.

\subsection{Asymptoticly valid confidence intervals via de-biasing}
\label{sec:debias}

We can leverage the results in Theorem~\ref{thm:posterior_contraction_mssl} and Corollary~\ref{cor:B_Psi_recovery} to derive asymptotically valid confidence intervals for individual parameters $\beta_{j,k}$ using de-biasing arguments similar to \citet{Rhyne2019} and \citet{Jankova2019}.
At a high-level, we first ``invert'' the Karush-Kuhn-Tucker conditions of the mSSL optimization problem to find the bias induced by the penalty.
Then we construct a new estimator that subtracts the bias from the mSSL MAP estimator. 
With an additional sparsity assumption on the design matrix $X,$ we can show that this ``de-biased'' estimator has an asymptotic normal distribution.

To construct the de-biased estimator, first denote the MAP estimates of $(B,\theta, \Omega,\eta)$ by $(\hat{B}, \hat{\theta}, \hat{\Omega}, \hat{\eta}).$
Fixing $\theta = \hat{\theta}$ and $\eta = \hat{\eta}$ at their MAP values, we know that
\begin{align}
\label{eq:debias_B_Omega_problem}
(\hat{B}, \hat{\Omega}) = \argmin_{B,\Omega}\left\{-\frac{n}{2}\log\left\lvert \Omega\right\rvert + \frac{1}{2}\tr\left[(Y - XB)^{\top}(Y-XB)\Omega\right] - \log\pi(B\vert\hat{\theta}) - \log\pi(\Omega \vert \hat{\eta})\right\}.
\end{align}
The KKT conditions for~\eqref{eq:debias_B_Omega_problem} tell us that
\begin{align}
-\hat{\Omega}^{-1}+(Y-X\hat{B})^\top(Y-X\hat{B})/n+\hat{z}&=0\label{eq:kkt_omega}\\
-X^\top(Y-X\hat{B})\hat{\Omega}+\hat{\kappa}&=0 \label{eq:kkt_B}
\end{align}
where $\hat{z}$ and $\hat{\kappa}$ are, respectively, sub-gradients of the penalties in~\eqref{eq:debias_B_Omega_problem} (i.e., the log-prior densities).

Now denote $\hat{\Phi} = X^{\top}X/n$ and suppose that we can find an approximate inverse $\hat{\Theta}$ such that $\lVert \hat{\Theta}\hat{\Phi} - I \rVert_{2} = o_{p}(\epsilon_{n}/\sqrt{n}).$ 
Then we can rewrite Equations~\eqref{eq:kkt_omega} and~\eqref{eq:kkt_B} as
\begin{align}
    \hat\Theta\hat\kappa\hat\Omega^{-1}+n\hat\Theta\hat\Phi(\hat B-B_0)&=\hat\Theta X^\top (Y-XB)\\
    (\hat B-B_0)+\hat\Theta \hat\kappa\hat\Omega^{-1}/n&=\hat\Theta X^\top (Y-XB)/n-\Delta/\sqrt n\label{eq:inverted_KKT_B}
\end{align}
where $\Delta=\sqrt{n}(\hat\Theta\hat\Phi-I)(\hat B-B_0)$.
Now, define the new ``de-biased'' estimator as
\begin{equation}
\label{eq:debiased_B_definition}
\hat B_b=\hat{B}+\hat\Theta\kappa \hat\Omega^{-1}/n=\hat B + \hat\Theta X^\top (Y-X\hat B)/n,
\end{equation}
which takes a similar form as the de-biased estimators studied by \citet{Rhyne2019}.
Theorem~\ref{thm:debias_ci} shows that $\hat{B}_{b}$ is asymptotically normal.
\begin{theorem}[Debiased confidence interval]
\label{thm:debias_ci}
Suppose $X$ has finite restricted eigenvalue and that one can obtain an approximate inverse of $\hat{\Phi} = (X^\top X)/n,$ denoted $\hat{\Theta}$ such that $\lVert \hat{\Theta}(X^{\top}X)/n - I\rVert_{2} = o_{p}(\epsilon_{n}/\sqrt{n}).$
Then,
\begin{align}
    \sqrt{n}(\hat B_b-B_0)_{\cdot,k}&\stackrel{\text{weakly}}{\rightarrow} \N_p(0,(\Omega^{-1})_{k,k}\hat\Theta\hat\Phi\hat{\Theta}^\top)
\end{align}
\end{theorem}

\begin{proof}
Denote the (true) residual as $E:=Y-XB$ and its column as $E_j$, we have $E\sim \N(0,\Omega^{-1})$ by assumption. 
Observe that
\begin{align}
    \begin{split}
        \sqrt{n}(\hat B_b-B_0)&=W-\Delta,\\
        W_k=\hat\Theta X^\top E_j/\sqrt n& \to \N_p(0,(\Omega^{-1})_{k,k}\hat\Theta\hat\Phi\hat{\Theta}^\top)\\
        \Delta &= \sqrt{n}(\hat\Theta(X^\top X)/n-I)(\hat B-B_0)
    \end{split}
\end{align}

The second line is due to CLT on the error term $E_j$. 
We only need to show $||\Delta||=o_p(1)$. 
Using the contraction result of mSSL (Corollary~\ref{cor:B_Psi_recovery}), we have $||\hat B-B_0||_F=o_p(\epsilon_n).$
By assumption on the approximated inverse $\hat\Theta$, we have $||\sqrt{n}(\hat\Theta(X^\top X)/n-I)(\hat B-B_0)||_2\le \sqrt{n}||\hat\Theta(X^\top X)/n-I||_2||\hat B-B_0||_F= o_p(1).$
\end{proof}

\begin{remark}[Approximate invertibility of $X^{\top}X/n$]
To form de-biased intervals for (univariate) LASSO, \citet{javanmard2018debiasing} \citep[also][for multivariate LASSO]{Rhyne2019} makes the same, arguably strong, assumption that $X^{\top}X/n$ be approximately invertible. 
In their procedure, they recommend estimating the approximate inverse using the graphical LASSO \citep[GLASSO;][]{Friedman2008}. 
In Section S2 of the Supplemental Materials, we show the spike-and-slab LASSO analog of GLASSO, which we call gSSL and which can be fit with a slightly modified mSSL procedure, contracts fast enough to use for estimating the approximate inverse $\hat{\Theta}$ \citep[see also][]{Gan2019_unequal}.

\end{remark}

It is tempting to follow a similar procedure to get a debiased estimator of $\Omega.$
Specifically, observe that Equation~\eqref{eq:kkt_omega} implies
\begin{align}
    -\hat{\Omega}+\hat\Omega\hat{z}\hat\Omega=-\hat\Omega(Y-X\hat{B})^\top(Y-X\hat{B})\hat\Omega/n.
\end{align}
Adding $2\hOmega-\Omega_0$ to each side, we see
\begin{align*}
    \hOmega+\hOmega\hat z\hOmega-\Omega_0=-\Omega_0 ((Y-XB)^\top(Y-XB)/n-\Omega_0^{-1})\Omega_0+\text{rem}_1+\text{rem}_2
\end{align*}
where
\begin{align*}
\text{rem}_1&=-(\hOmega-\Omega_0)((Y-X\hat{B})^\top(Y-X\hat{B})/n-\Omega_0^{-1})\Omega_0-(\hOmega(Y-X\hat{B})^\top(Y-X\hat{B})-I)(\hOmega-\Omega_0)\\
\text{rem}_2&=-\Omega_0((Y-X\hat{B})^\top(Y-X\hat{B})-(Y-XB_0)^\top(Y-XB_0))/n\Omega_0    
\end{align*}

If we could show that $\ell_{\infty}$ norm of the two remainder terms were $o_{p}(1/\sqrt{n}),$ then the $(k,k')$ entry of $\sqrt{n}( (Y-XB)^{\top}(Y - XB)/n - \Omega_{0}^{-1})$ would converge weakly to $\N(0,\sigma^{2}_{k,k'}),$ where $\sigma^{2}_{k,k'} = \Omega_{0,k,k}\Omega_{0,k,k'}+\Omega_{k,k'}^2$. 
We could further form a ``de-biased'' estimator of $\Omega$ as
\begin{equation}
\label{eq:debiased_Omega_definition}
 \hat\Omega_b=2\hat\Omega-\hat\Omega(Y-X \hat B)^\top(Y-X \hat B)\hat\Omega/n.
\end{equation}
Unfortunately, establishing an analog of Theorem~\ref{thm:debias_ci} for the estimator in~\eqref{eq:debiased_Omega_definition} requires a posterior contraction rate in $\ell_{\infty}-$norm  that is faster than that implied by our results for the Frobenius norm.
Though we cannot establish the asymptotic result, in Section~\ref{sec:simulation}, we nevertheless assess the finite-sample performance of intervals based on the de-biased estimator $\hat{\Omega}_{b}$ and the plug-in estimator of the anticipated asymptotic variance $\hat{\sigma}_{k,k'} = \hat{\Omega}_{k,k}\hat{\Omega}_{k,k'} + \hat{\Omega}_{k,k'}^{2}.$

\section{Finite sample uncertainty quantification\label{sec:simulation}}
We performed a simulation study to compare the frequentist operating characteristics of several approaches for constructing mSSL uncertainty intervals.
We specifically compared forming marginal 95\% confidence intervals based on debiasing (i.e., using Theorem~\ref{thm:debias_ci}) to three variants of the Bayesian bootstrap.
For each variant, we first run dynamic posterior exploration to obtain point estimates $\hat{B}, \hat{\theta}, \hat{\Omega}$ and $\hat{\eta}.$
Then, we repeatedly solve the the randomized optimization problem in Equation~\eqref{eq:BB_main_equation} to obtain 500 bootstrap re-samples of $(B, \Omega).$
The first Bayesian bootstrap variant, which we call \texttt{WLB} after \citet{Newton1994}'s ``weighted likelihood bootstrap,'' sets $w_{0} = 0$ and $B_{0} = \mathbf{0}_{p \times q}$ and $\Omega_{0} = \mathbf{0}_{q \times q}.$
The second variant, which we call \texttt{WBB} after \citet{Newton1994}'s ``weighted Bayesian bootstrap,'' also sets $B_{0}$ and $\Omega_{0}$ to the appropriately-sized zero matrices but instead draws $w_{0} \sim \textrm{Gamma}(1,1).$
The final variant, which we call \texttt{BB-jitter}, draws $w_{0} \sim \textrm{Gamma}(1,1)$, $B_{0,j,k} \sim \textrm{Laplace}(\lambda_0)$ and $\Omega_{0} \sim \textrm{Laplace}(\xi_0)$.
\texttt{BB-jitter} is an analog to the procedure introduced in \citet{Nie2020} that randomly re-centers or ``jitters'' the log-prior in Equation~\eqref{eq:BB_main_equation}.

\subsection{Simulation design}

\begin{figure}[h]
\centering
\begin{subfigure}[b]{0.19\textwidth}
\centering
\includegraphics[width = \textwidth]{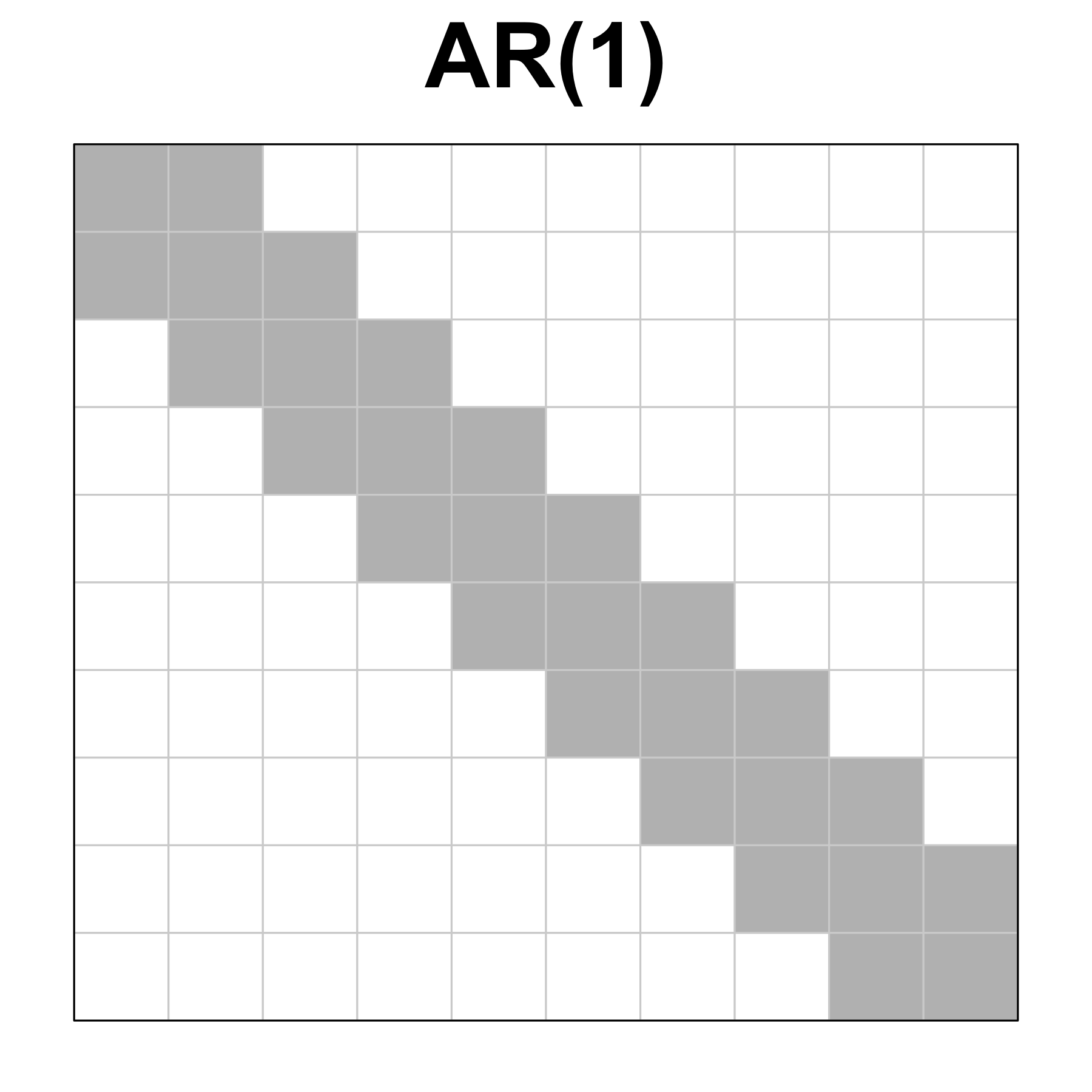}
\end{subfigure}
\begin{subfigure}[b]{0.19\textwidth}
\centering
\includegraphics[width = \textwidth]{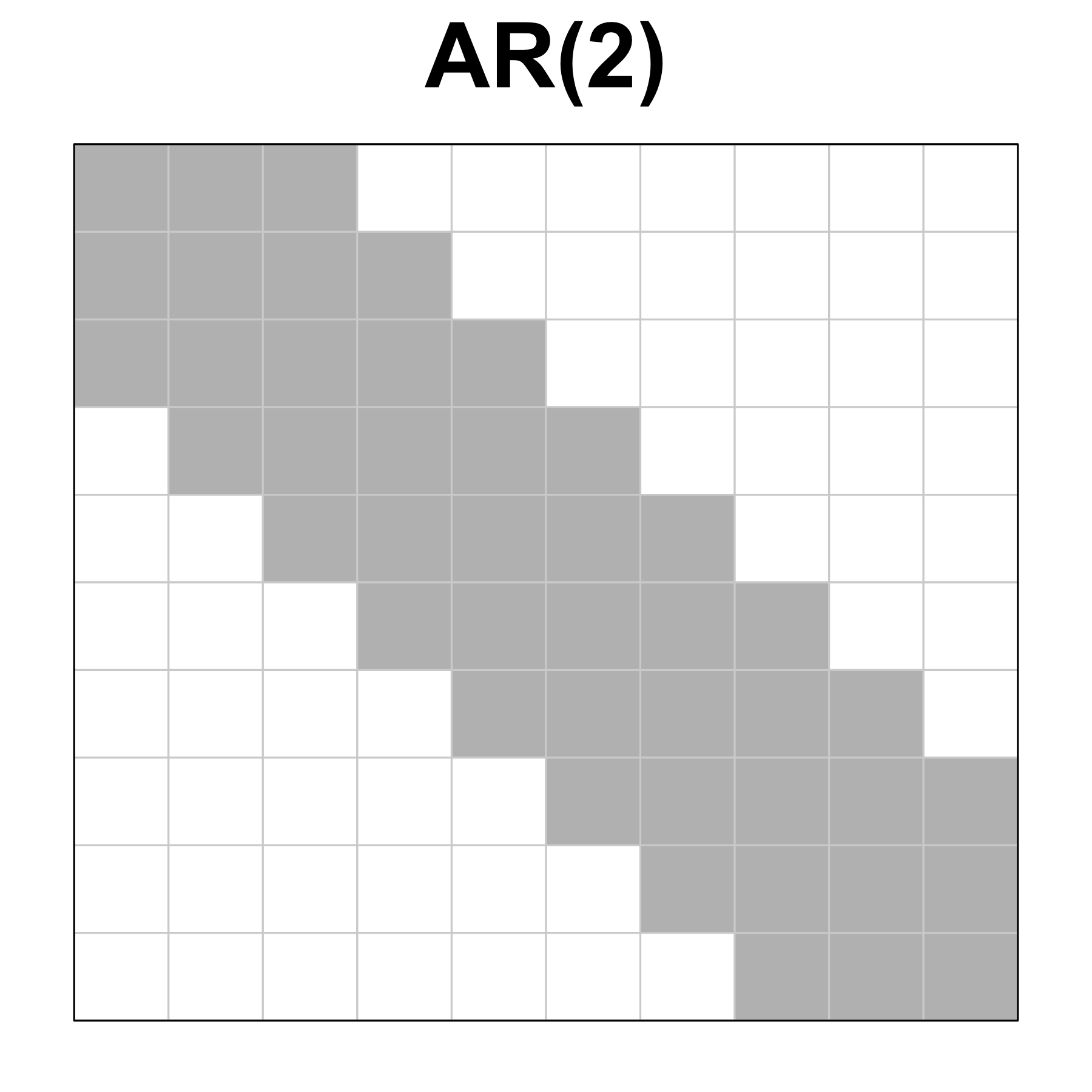}
\end{subfigure}
\begin{subfigure}[b]{0.19\textwidth}
\centering
\includegraphics[width = \textwidth]{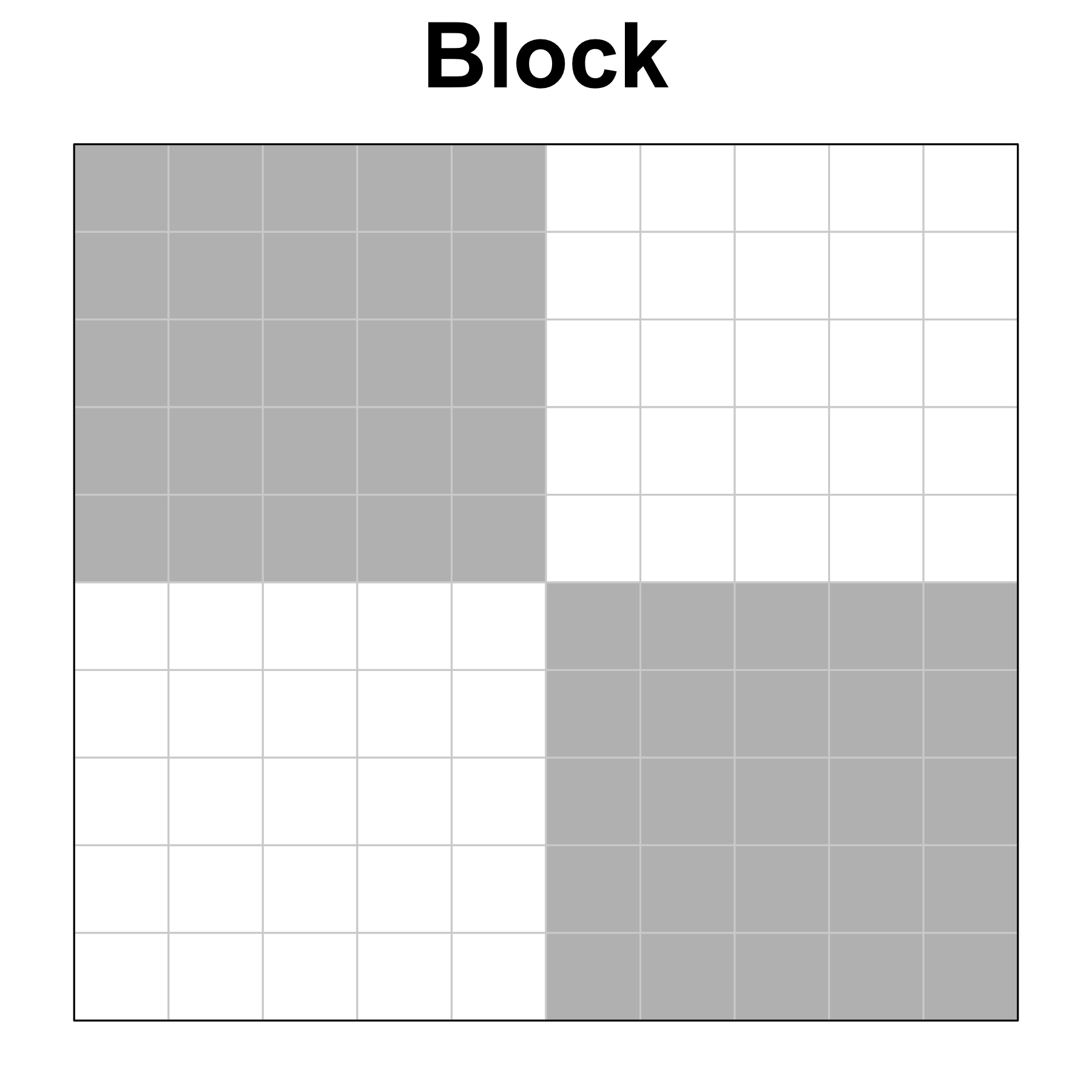}
\end{subfigure}
\begin{subfigure}[b]{0.19\textwidth}
\centering
\includegraphics[width = \textwidth]{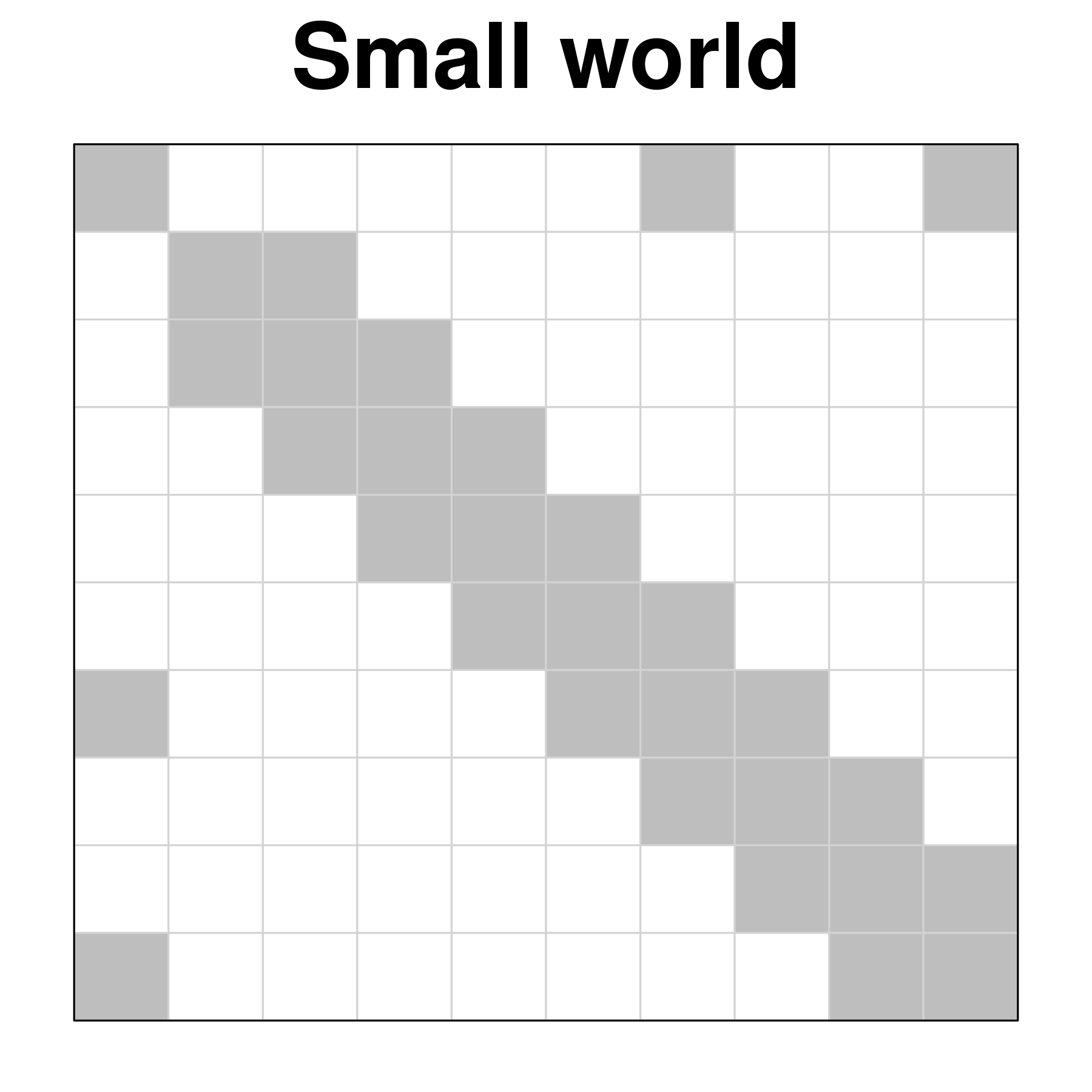}
\end{subfigure}	
\begin{subfigure}[b]{0.19\textwidth}
\centering
\includegraphics[width = \textwidth]{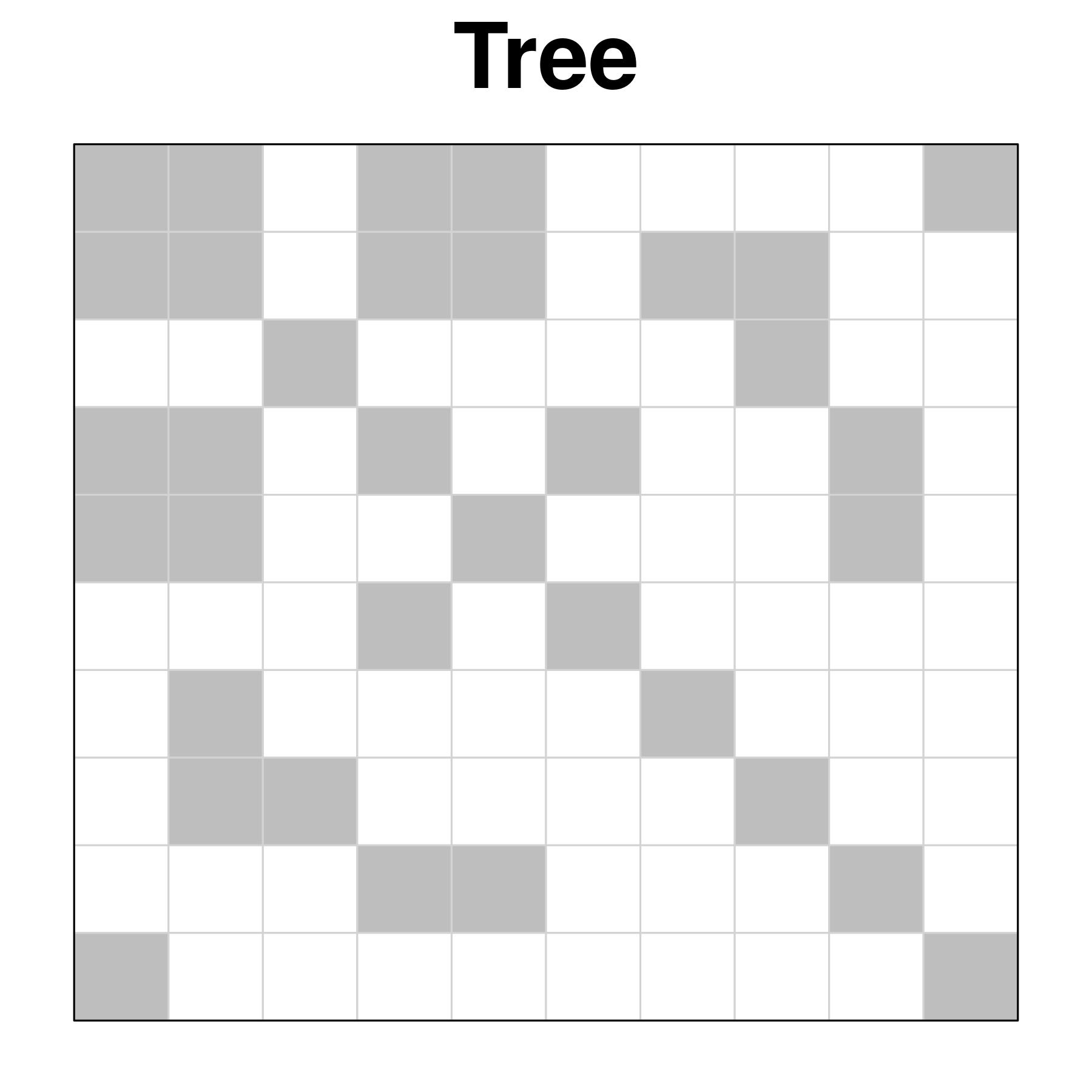}
\end{subfigure}	

\begin{subfigure}[b]{0.19\textwidth}
\centering
\includegraphics[width = \textwidth]{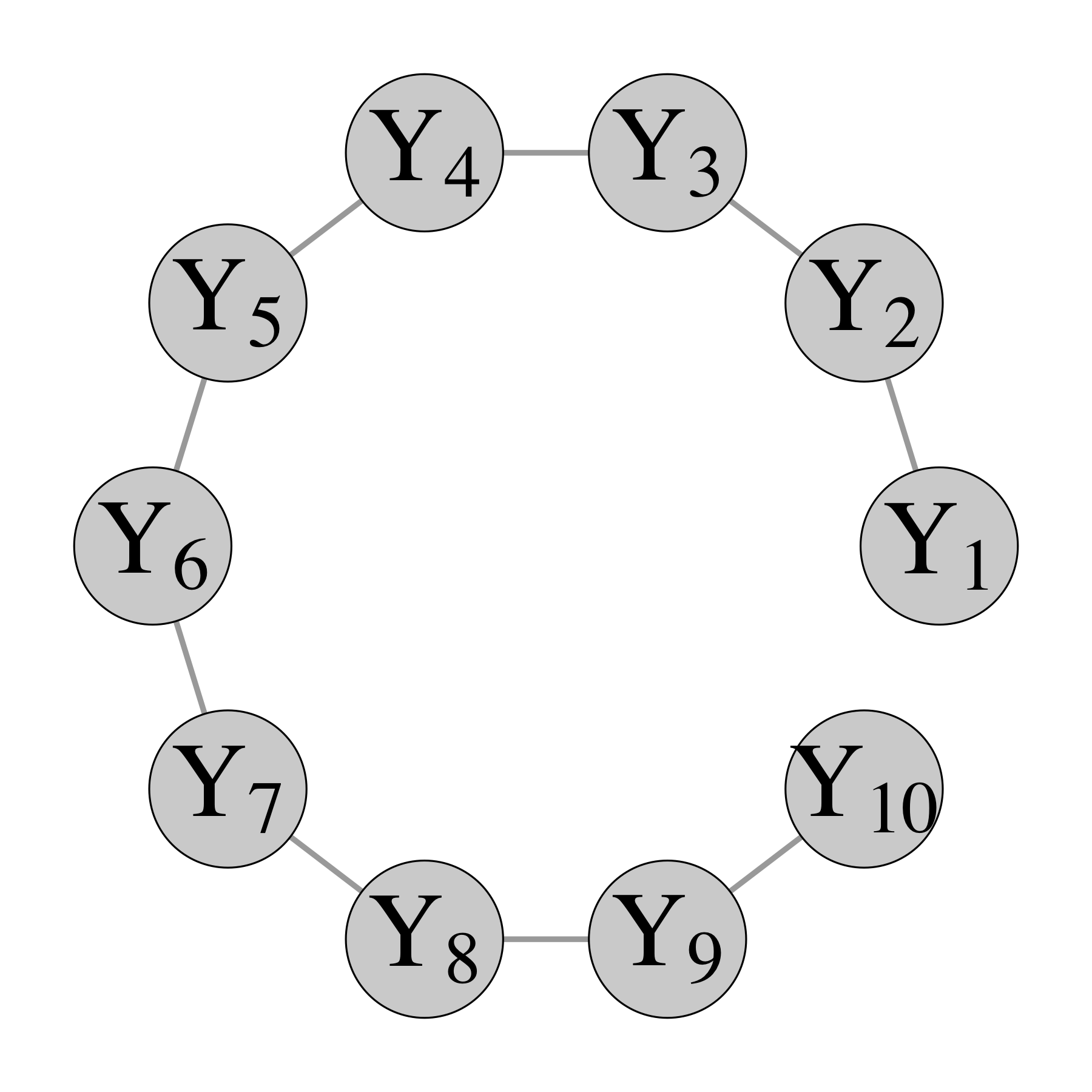}
\end{subfigure}
\begin{subfigure}[b]{0.19\textwidth}
\centering
\includegraphics[width = \textwidth]{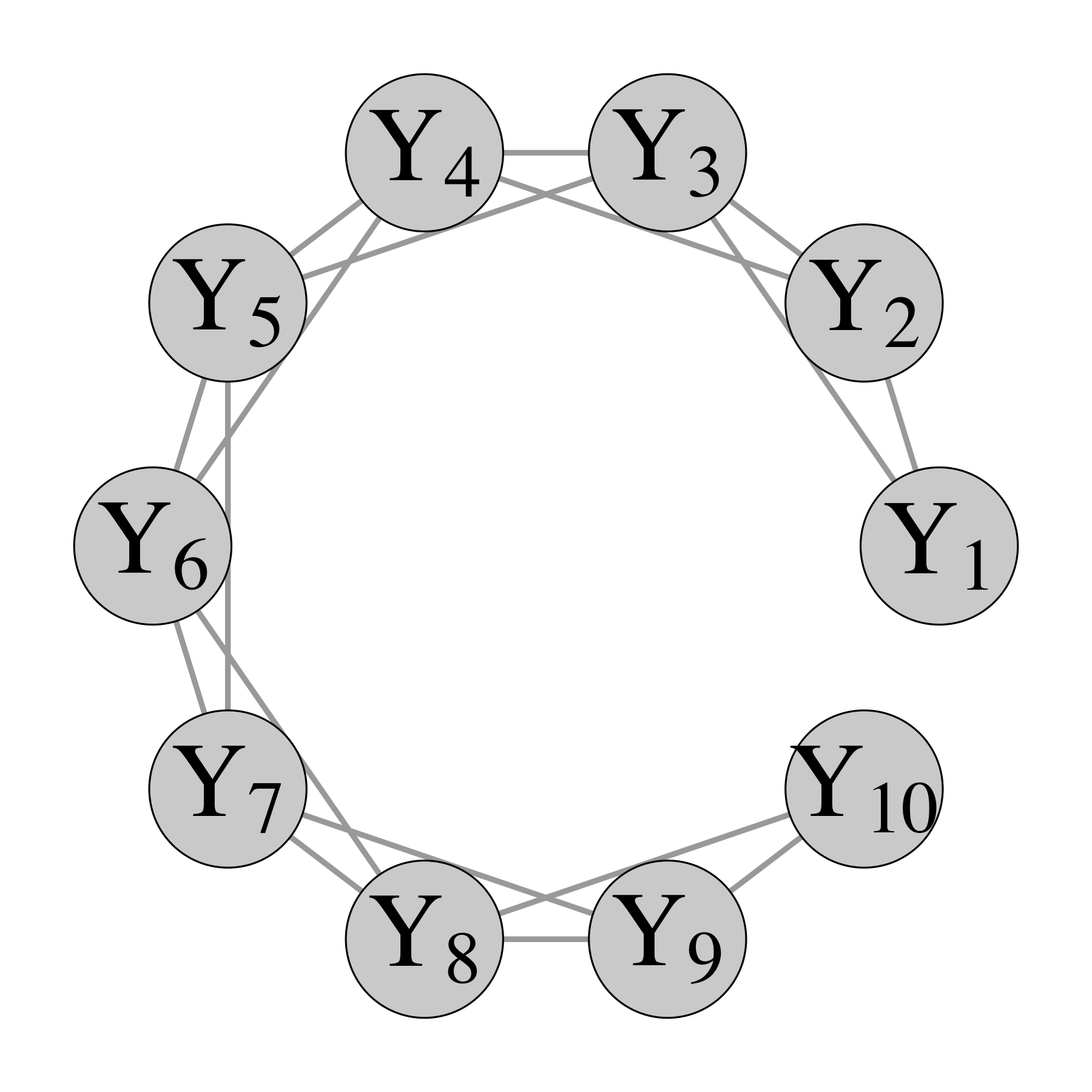}
\end{subfigure}
\begin{subfigure}[b]{0.19\textwidth}
\centering
\includegraphics[width = \textwidth]{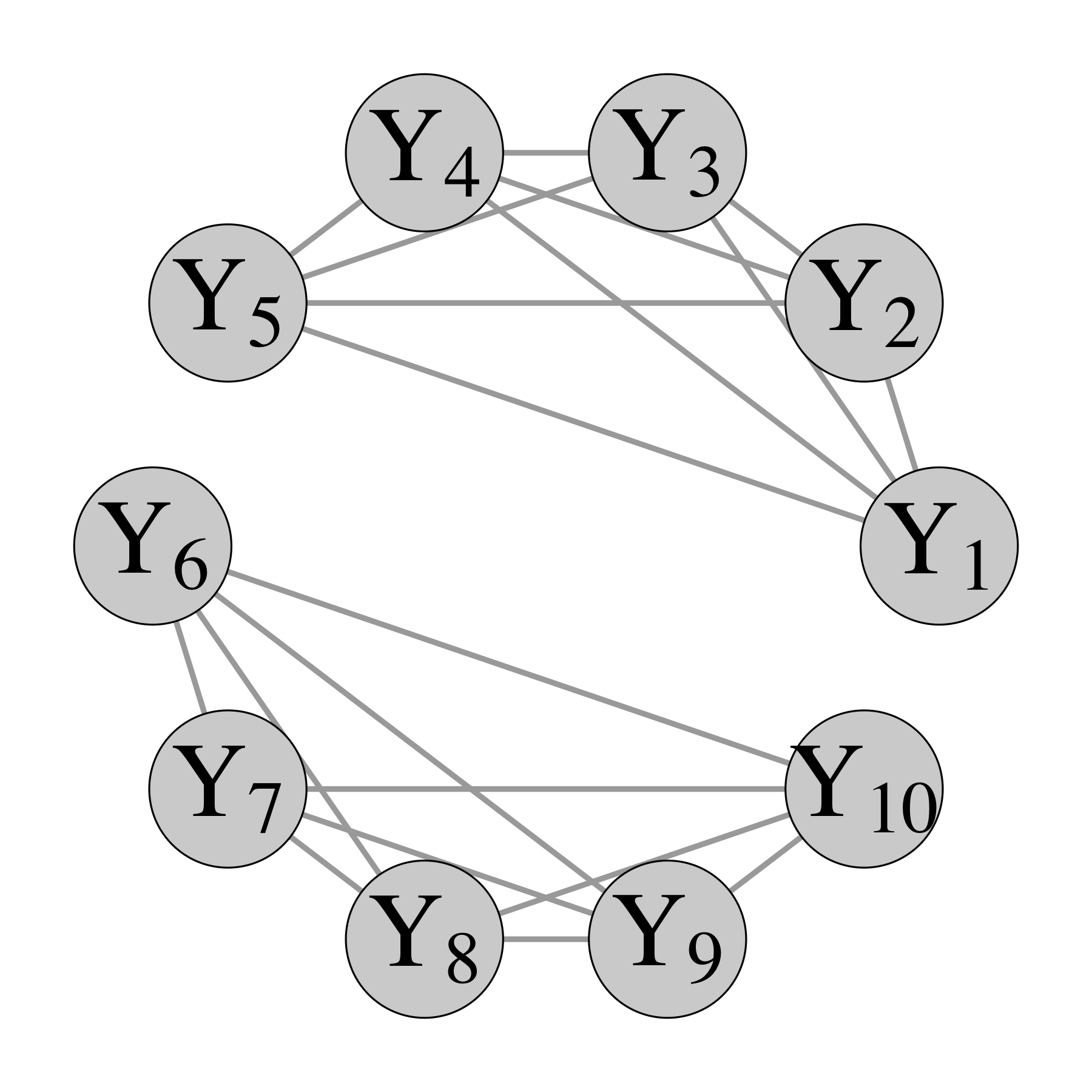}
\end{subfigure}
\begin{subfigure}[b]{0.19\textwidth}
\centering
\includegraphics[width = \textwidth]{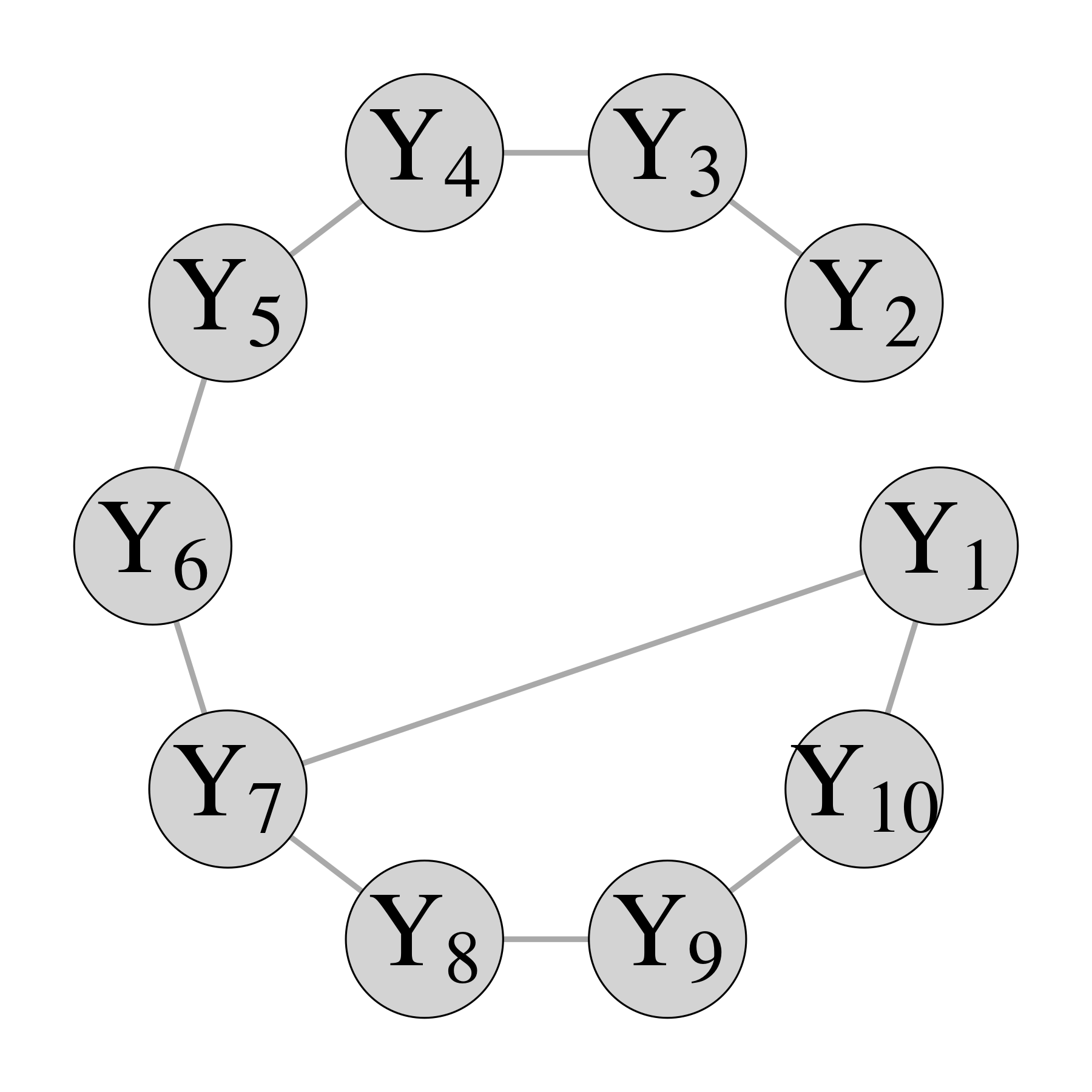}
\end{subfigure}	
\begin{subfigure}[b]{0.19\textwidth}
\centering
\includegraphics[width = \textwidth]{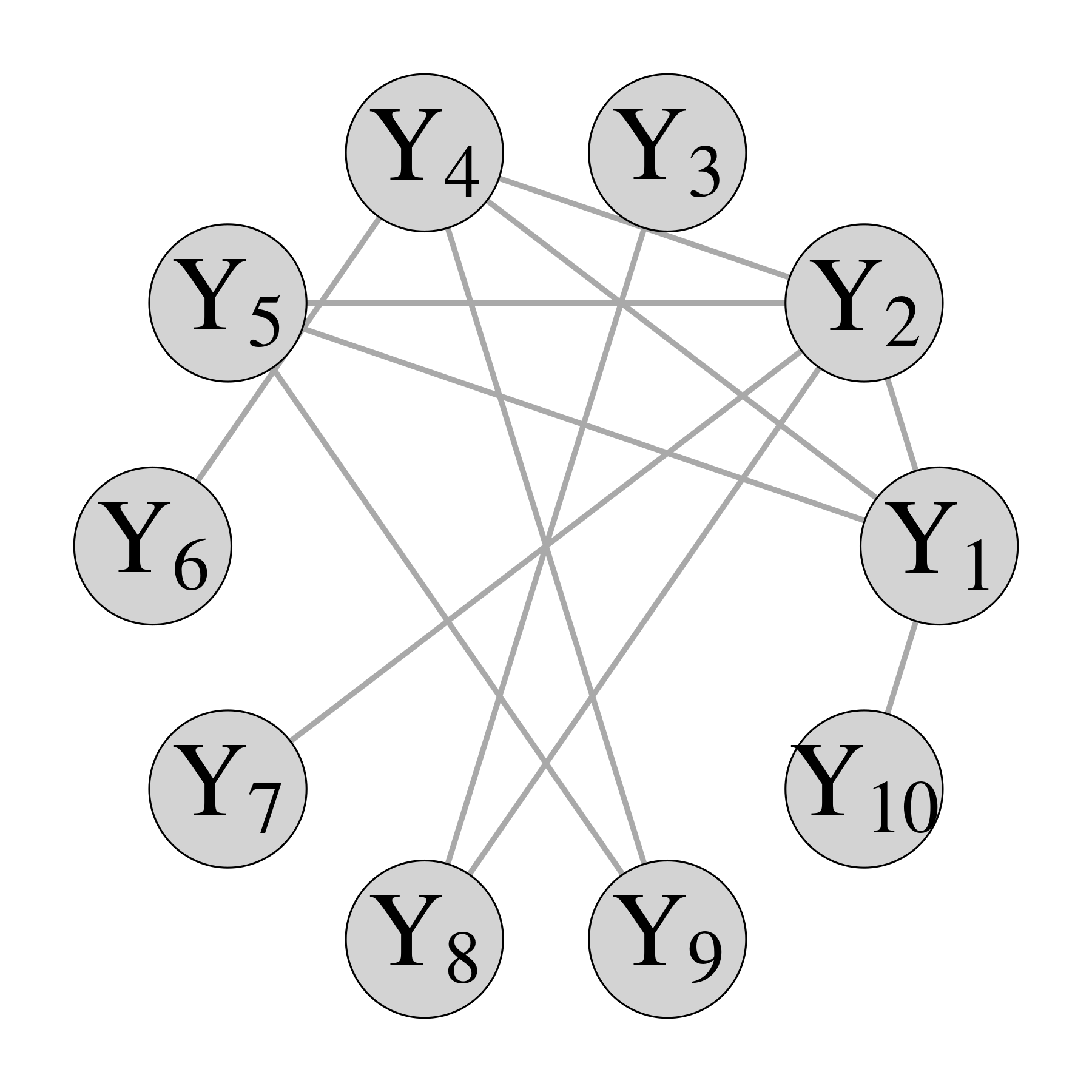}
\end{subfigure}	
\caption{Visualization of the supports of $\Omega$ for $q = 10$ under the five precision structures in simulation studies (top) and their graph representation (bottom). Gray cells indicate non-zero entries in $\Omega$ while white cells indicate zeros}
    \label{fig:simulation_design}
\end{figure}

We performed a total of 15 simulation studies, one for each combination of three problem dimensions and seven sparsity patterns in $\Omega.$ 
We considered three settings of problem dimension $(n,p,q) = (100, 10, 10), (100, 20, 30)$, and $(400, 100, 30).$ 
Figure~\ref{fig:simulation_design} shows cartoon illustrations of the different sparsity patterns for $\Omega$; we defer precise formulas to Section S3.1 of the Supplementary Materials.
For each choice of dimension and $\Omega,$ we drew a random $p \times q$ design matrix $X$ with independent standard normal entries and a random, sparse $p \times q$ matrix $B$ with $pq/5$ non-zero entries drawn uniformly from the interval $[-2,2].$
We generated 100 synthetic datasets from the model in Equation~\eqref{eq:general_model} for each combination of $X, B,$ and $\Omega.$

For each simulated dataset, we formed 95\% uncertainty intervals for each parameter $\beta_{j,k}$ and $\omega_{k,k'}.$
We then computed the proportion of intervals containing the true data-generating parameter values.
By averaging these proportions across the different simulation replications, we estimated the marginal frequentist coverage of our intervals.
Beyond coverage, we also computed the average ratio of uncertainty interval lengths relative to \texttt{debiasing}.
We finally computed the average interval score (across all entries in $B$ and across all entries in $\Omega$).
Note that the interval score is a proper scoring rule that balances coverage and interval length \citep{gneiting2007strictly}. 

\citet{Deshpande2019} performed variable and covariance selection by identifying the non-zero entries in the MAP estimates of $B$ and $\Omega.$
Given marginal uncertainty intervals for each $\beta_{j,k}$ and $\omega_{k,k'},$ a natural alternative involves identifying those intervals not containing zero.
We compared the sensitivity, specificity, precision, and overall $F_{1}$ score for such selection of such interval-based selection to \citet{Deshpande2019}'s original estimate-based selection.

\subsection{Simulation results}

Table~\ref{tab:simulation_on_debias_b} and \ref{tab:simulation_on_debias_omega} shows the variable selection and inference performance of method we tested when $(n,p,q) = (400, 100, 30)$. 
Since the intervals formed by \texttt{debiasing} and each version of the Bayesian bootstrap are based on the \texttt{mSSL} point estimate, we see that the former methods inherit the latter's excellent recovery of $B$'s support.
The same is not true for $\Omega,$ however: interval-based selection appears to be better than estimate-based selection in some settings and to be worse in others.
Both de-biasing and the Bayesian bootstraps produced uncertainty intervals with higher-than-nominal frequentist coverage for both $B$ and $\Omega.$ with the former producing much longer intervals than the latter.
Generally speaking \texttt{WLB} produced the shortest intervals followed by \texttt{WBB}, \texttt{BB-jitter}, and then \texttt{debiasing}.
On this basis, we would recommend using \texttt{WLB} to quantify the sampling uncertainty of the MAP estimate returned by \texttt{mSSL}.
\begin{table}[H]
    \centering
    \caption{Support recovery and uncertainty interval coverage for $B$ when $(n,p,q) = (400, 100, 30)$ average across 100 simulated datasets. Best performance is bold-faced.}
    \label{tab:simulation_on_debias_b}
    \scriptsize
    \begin{tabular}{lccccccc}
        \hline
        Method & SEN & SPEC & PREC & F1 & Relative length& Coverage & Interval score  \\ \hline
        \multicolumn{8}{c}{$AR(1)$ model} \\ \hline
        \texttt{debiasing} & \textbf{0.94}  & \textbf{1}  & \textbf{1}  & \textbf{0.97}  & 1 & 0.95  & 3.9 \\
        \texttt{WLB} & \textbf{0.94}  & \textbf{1}  & \textbf{1} & \textbf{0.97}  & \textbf{0.12}  & 0.97  & \textbf{0.52}  \\
        \texttt{WBB} & \textbf{0.94}  & \textbf{1}  & \textbf{1}  & \textbf{0.97}  & 0.25  & \textbf{0.98}  & 1.4  \\
        \texttt{BB-jitter} & 0.92  & \textbf{1}  & \textbf{1}  & 0.96  & 0.29  & \textbf{0.98}  & 1.1  \\
        \texttt{mSSL} & \textbf{0.94} & \textbf{1}  & \textbf{1}  & \textbf{0.97}  & NA & NA & NA \\
        \hline
        \multicolumn{8}{c}{$AR(2)$ model} \\ \hline
        \texttt{debiasing} & \textbf{0.90} & \textbf{1} & \textbf{1} & \textbf{0.94} & 1 & 0.95 & 5.2 \\
        \texttt{WLB} & \textbf{0.90} & \textbf{1} & \textbf{1} & \textbf{0.94} & \textbf{0.15} & \textbf{0.97} & \textbf{0.91} \\
        \texttt{WBB} & \textbf{0.90} & \textbf{1} & \textbf{1} & \textbf{0.94} & 0.25 & \textbf{0.97} & 1.8 \\
        \texttt{BB-jitter} & 0.85 & \textbf{1} & \textbf{1} & 0.92 & 0.38 & \textbf{0.97} & 1.5\\
        \texttt{mSSL} & \textbf{0.90}  & \textbf{1} & \textbf{1} & \textbf{0.94} & NA & NA & NA \\
        \hline
        \multicolumn{8}{c}{Block model} \\ \hline
        \texttt{debiasing} & \textbf{0.93} & \textbf{1} & \textbf{1} & \textbf{0.96} & 1 & 0.95 & 3.9 \\
        \texttt{WLB} & \textbf{0.93} & \textbf{1} & \textbf{1} & \textbf{0.96} & \textbf{0.14} & 0.97 & \textbf{0.61} \\
        \texttt{WBB} & \textbf{0.93} & \textbf{1} & \textbf{1} & \textbf{0.96} & 0.28 & \textbf{0.98} & 1.6 \\
        \texttt{BB-jitter} & 0.90 & \textbf{1} & \textbf{1} & 0.95 & 0.25 & \textbf{0.98} & 1.2 \\
        \texttt{mSSL} & \textbf{0.93}  & \textbf{1} & \textbf{1} & \textbf{0.96} & NA & NA & NA\\
        \hline
        \multicolumn{8}{c}{Small world} \\ \hline
        \texttt{debiasing} & \textbf{0.94} & \textbf{1} & \textbf{1} & \textbf{0.97} & 1 & 0.95 & 3.1 \\
        \texttt{WLB} &\textbf{0.94} & \textbf{1} & \textbf{1} & \textbf{0.97} & \textbf{0.15} & \textbf{0.98} & \textbf{0.47} \\
        \texttt{WBB} &\textbf{0.94} & \textbf{1} & \textbf{1} & \textbf{0.97} & 0.25 & \textbf{0.98} & 1.3\\
        \texttt{BB-jitter} &0.92 & \textbf{1} & \textbf{1} & 0.96 & 0.33 & \textbf{0.98} & 1\\
        \texttt{mSSL} &\textbf{0.94} & \textbf{1} & \textbf{1} & \textbf{0.97} &NA&NA& NA \\
        \hline
        \multicolumn{8}{c}{Tree model} \\ \hline
        \texttt{debiasing} & \textbf{0.93} & \textbf{1} & \textbf{1} & \textbf{0.96} & 1 & 0.95 & 3.9 \\
        \texttt{WLB} &\textbf{0.93} & \textbf{1} & \textbf{1} & \textbf{0.96} & \textbf{0.13} & 0.97 & \textbf{0.55} \\
        \texttt{WBB} & \textbf{0.93} & \textbf{1} & \textbf{1} & \textbf{0.96} & 0.35 & \textbf{0.98} & 1.4\\
        \texttt{BB-jitter} &0.91 & \textbf{1} & \textbf{1} & 0.95 & 0.37 & \textbf{0.98} & 1.1 \\
        \texttt{mSSL} &\textbf{0.93} & \textbf{1} & \textbf{1} & \textbf{0.96} & NA & NA & NA \\
        \hline
    \end{tabular}
\end{table}

\begin{table}[H]
    \centering
       \caption{Support recovery and uncertainty interval coverage for $\Omega$ when $(n,p,q) = (400, 100, 30)$ average across 100 simulated datasets. Best performance is bold-faced.}
    \label{tab:simulation_on_debias_omega}
    \scriptsize
    \begin{tabular}{lccccccc}
        \hline
        Method & SEN & SPEC & PREC & F1 & Length& Coverage & Interval score  \\ \hline
        \multicolumn{8}{c}{$AR(1)$ model} \\ \hline
        \texttt{debiasing} & \textbf{1} & \textbf{1} & 0.99 & \textbf{1} & 1 & 0.95 & 11 \\
        \texttt{WLB} & \textbf{1} & \textbf{1} & 0.99 & 0.99 & \textbf{0.07} & \textbf{1} & \textbf{0.86} \\
        \texttt{WBB} & \textbf{1} & \textbf{1} & 0.99 & 0.99 & 0.15 & \textbf{1} & 1.3 \\
        \texttt{BB-jit} & \textbf{1} & \textbf{1} & \textbf{1} & \textbf{1} & 0.48 & \textbf{1} & 1.4 \\
        \texttt{mSSL} & \textbf{1} & \textbf{1} & 0.99 & 0.99  & NA & NA & NA\\
        \hline
        \multicolumn{8}{c}{$AR(2)$ model} \\ \hline
        \texttt{debiasing} & \textbf{1} & \textbf{1} & \textbf{0.99} & \textbf{1} & 1 & 0.95 & 4.2\\
        \texttt{WLB} & \textbf{1} & \textbf{1} & \textbf{0.99} & \textbf{1} & \textbf{0.12} & 0.99 & \textbf{0.66} \\
        \texttt{WBB} & \textbf{1} & \textbf{1} & \textbf{0.99} & \textbf{1} & 0.33 & \textbf{1} & 2.2 \\
        \texttt{BB-jit} & \textbf{1} & \textbf{1} & \textbf{0.99} & \textbf{1} & 0.17 & \textbf{1} & 1.9 \\
        \texttt{mSSL} & \textbf{1} & \textbf{1} & \textbf{0.99}  & \textbf{1}  & NA & NA & NA \\
        \hline
        \multicolumn{8}{c}{Block model} \\ \hline
        \texttt{debiasing} & 0.26 & \textbf{1}& 0.99 & 0.41 & 1 & \textbf{0.93} & 8.1 \\
        \texttt{WLB} & \textbf{0.40} & \textbf{1}& 0.99 & \textbf{0.57} & \textbf{0.20} & 0.61 & \textbf{3.5} \\
        \texttt{WBB} & \textbf{0.40} & \textbf{1}& 0.99 & \textbf{0.57} & 0.35 & 0.73 & 3.8 \\
        \texttt{BB-jit} & 0.32 & \textbf{1}& \textbf{1}& 0.49 & 0.71 & 0.77 & 3.8 \\
        \texttt{mSSL} & \textbf{0.40}  & \textbf{1} & 0.99 & \textbf{0.57}  & NA & NA & NA\\
        \hline
        \multicolumn{8}{c}{Small world} \\ \hline
        \texttt{debiasing} & 0.64 & \textbf{1}& 0.94 & 0.76 & 1 & 0.95 & 19 \\
        \texttt{WLB} &0.65 & 0.99 & 0.82 & 0.73 & \textbf{0.07} & 0.90 & \textbf{2.9} \\
        \texttt{WBB} &0.65 & 0.99 & 0.82 & 0.73 & 0.65 & 0.91 & 3.7\\
        \texttt{BB-jit} &0.64 & 0.99 & 0.89 & 0.74 & 0.51 & \textbf{0.98} & 3.9\\   
        \texttt{mSSL} &\textbf{0.90} & \textbf{1}& \textbf{0.97} & \textbf{0.93} & NA & NA & NA \\
        \hline
        \multicolumn{8}{c}{Tree model} \\ \hline
        \texttt{debiasing} & 0.88 & \textbf{1}& \textbf{0.99} & 0.93 & 1 & 0.96 & 15 \\
        \texttt{WLB} &\textbf{0.90} & \textbf{1}& 0.97 & 0.93 & \textbf{0.071} & 0.96 & \textbf{1.4} \\
        \texttt{WBB} &\textbf{0.90} & \textbf{1}& 0.97 & 0.93 & 0.16 & 0.97 & 2.9\\
        \texttt{BB-jit} &0.89 & \textbf{1}& \textbf{0.99} & \textbf{0.94} & 0.19 & \textbf{0.99} & 3.0\\
        \texttt{mSSL} &\textbf{0.90} & \textbf{1}& 0.97 & 0.93 & NA & NA & NA\\
        \hline
    \end{tabular}
\end{table}

\nocite{Roverato2002, Lenkoski2013,WattsStrogatz1998}

\section{Discussion\label{sec:discussion}}
Our theoretical results show that full posterior distribution of a slightly modified version of \citet{Deshpande2019}'s mSSL procedure concentrates around the true data generating parameter. 
Our analysis, with minimal modifications, provides the posterior contraction rate for the so-called ``graphical spike-and-slab LASSO'' in which one simply wishes to estimate a sparse precision matrix using spike-and-slab LASSO priors on the off-diagonal elements (Section S1.4 of the Supplemental Materials).
Beyond entry-wise sparsity, it is not difficult to extend \citet{Bai2020_groupSSL}'s group SSL prior to the multiple outcome setting.
Such an extension would facilitate, for instance, sparse additive modeling of multiple correlated outcomes. 
From a computational perspective, we can deploy a similar ECM algorithm as the mSSL.
Our proof strategy for the mSSL can also be modified to establish the posterior contraction rate for the group version of the mSSL.
Specifically our key packing number lemma can be used for a ``effectively group sparse'' set that has the form of $\mathcal{A}\times \{a:||a||_2<\delta,a\in \R^r\}.$

It would also be interesting to study the theoretical properties of \citet{Moran2021_biclustering}'s SSL biclustering posterior.
Their work involves a factor model with a multivariate Gaussian error structure that induces a similar KL-geometry as the mSSL model we studied here. 
Such theoretical analysis is complicated by the fact that \citet{Moran2021_biclustering}'s factor model is identifiable only up to rotation.



\bibliographystyle{apalike}
\bibliography{references}

\renewcommand{\theequation}{S\arabic{equation}}
\renewcommand{\thesection}{S\arabic{section}}  
\renewcommand{\thefigure}{S\arabic{figure}}  
\renewcommand{\thetable}{S\arabic{table}} 
\renewcommand{\themyLemma}{S\arabic{myLemma}} 
\renewcommand{\themyTheorem}{S\arabic{myTheorem}} 
\renewcommand{\themyProposition}{S\arabic{myProposition}}
\setcounter{equation}{0}
\setcounter{section}{0}
\setcounter{subsection}{0}
\setcounter{subsubsection}{0}
\setcounter{myLemma}{0}
\setcounter{myTheorem}{0}

\newpage
\begin{center}
	{\LARGE \textbf{Supplementary Materials}}
\end{center}

In Section~\ref{appendix:gSSL}, we consider Gaussian graphical modeling with spike-and-slab LASSO priors and establish concentration rates for the so-called ``gSSL'' posterior.
As noted in Remark 1 in the main text, we used these results to construct de-biased intervals for $\beta_{j,k}.$
We conclude in Section~\ref{appendix:experiments} with additional experimental details and results.


\section{Gaussian graphical modeling with the spike-and-slab LASSO}
\label{appendix:gSSL}
In this section, we consider a special case of our multi-outcome regression model in which there are no covariates (i.e., $p = 0$).
Specifically we model 
$$
\by \vert \Omega \sim \N(0, \Omega^{-1}).
$$
The gSSL procedure works by specifying spike-and-slab LASSO priors on the off-diagonal elements of $\Omega.$
We can obtain the MAP estimate of $\Omega$ using an EM algorithm that arises as a special case of the mSSL ECM algorithm that fixes $B = \bm{0}$ \citep[see also ][]{Gan2019_unequal}
We can further show that the so-called gSSL posterior concentrates using virtually the same strategy as in Section~\ref{appendix:proofs}.


\begin{theorem}[Contraction in log affinity of gSSL]
    \label{thm:log-affinity-gssl}
    For Gaussian graphical model $Y_i\sim \N(0,\Omega^{-1})$, with prior on $\Omega$ following the SSL prior and Assumptions A1, A5, and 
    \begin{itemize}
        \item[A3']: $\log(n)\lesssim \log(q)$ and $\max\{q,s_0^\Omega\}\log q/n\to 0$
    \end{itemize}

    There is some constant $M>0$ that does not depend on $n$ such that
    \begin{align}
        \sup_{\Omega\in\mathcal{H}_0}\E_0 \Pi\left(\Omega:\frac{1}{n}\sum \rho(f_i,f_{0i})\ge M\epsilon_n^{'2}|Y_1,\dots,Y_n\right) &\longrightarrow 0
    \end{align}
    where $\epsilon_n'=\sqrt{\max\{q,s_0^\Omega\}\log(q)/n},$ and $s_{0}^{\Omega}$ are the numbers of non-zero entries in off-digaonal entries in $\Omega$ and the log-affinity is defined as
    \begin{align*}
    \frac{1}{n}\sum \rho(f_i,f_{0,i}):=&-\log\left(\frac{|\Omega^{-1}|^{1/4}|\Omega_0^{-1}|^{1/4}}{|(\Omega^{-1}+\Omega_0^{-1})/2|^{1/2}}\right)
\end{align*}
\end{theorem}

The proof follows the proof of contraction of mSSL, by skipping all parts related to $B$ and establishing same technical lemmas for $\Omega$ only as follows. 

\begin{lemma}[KL condition for gSSL]
    Let $\epsilon'_n=\sqrt{\max\{q,s_0^\Omega\}\log(q)/n}.$ Then for all true parameters $\Omega_{0}$ we have
    \begin{align*}
        -\log\Pi\left[(\Omega):K(f_0,f)\le n\epsilon_n^{'2}, V(f_0,f)\le n\epsilon_n^{'2}\right]\le C_1n\epsilon_n^{'2},
    \end{align*}
    where $K$ and $V$ are, respectively, the Kullback-Leibler divergence

    \begin{equation}
        \label{eqn:KL_formula_gssl1}
        K(f_{0}, f) = \frac{n}{2}\left(\log\left(\frac{|\Omega_0|}{|\Omega|}\right)-q+\tr(\Omega_0^{-1}\Omega)\right).
        \end{equation}
        The KL variance between the same models is 
        \begin{equation}
        \label{eqn:KL_formula_gssl2}
        V(f_{0}, f) = \frac{n}{2}\left(\tr((\Omega_0^{-1}\Omega)^2)-2\tr (\Omega_0^{-1}\Omega) + q \right).
    \end{equation}

\end{lemma}

\begin{proof}

    Denote $\Omega-\Omega_0=\Delta_\Omega$, we consider the same event as in~\eqref{eqn:mssl_A1_due_to_Omega}

    \begin{equation*}
        \mathcal{A}_1=\left\{ \Omega: \tr((\Omega_0^{-1}\Omega)^2)-2\tr (\Omega_0^{-1}\Omega) + q\le \epsilon_n^2, \log\left(\frac{|\Omega_0|}{|\Omega|}\right)-q+\tr(\Omega_0^{-1}\Omega)\le \epsilon_n^{'2} \right\}
    \end{equation*}

    Which implies the KL event. By bounds in~\eqref{eq:B_event_for_Omega} we can get the results in the lemma. 
\end{proof}

\begin{lemma}[Dimension recovery of gSSL]
    \label{lemma:dimension_recovery_gssl}
    For a sufficiently large number $C_3'>0$, we have:
    \begin{align}
        \sup_{\Omega\in\mathcal{H}_0}\E_0\Pi\left(\Omega:|\nu_{\omega}(\Omega)|>C_3's^\star|Y_1,\dots,Y_n\right)&\to 0
    \end{align}
    where $s^\star=\max\{q,s_0^\Omega\}.$ 
\end{lemma}

\begin{proof}
    Same as proof in~\eqref{eq:Omega_dimension}, with modifications to $s^\star$. 
\end{proof}

\begin{lemma}[Sieve]
    \label{lemma:sieve_gssl}
There exist a sequence of sieve $\mathcal{F}_n$ such that it receive enough prior mass, i.e. $\Pi(\mathcal{F}_n^c)\le \exp(-C_2n\epsilon_n^{'2})$.
\end{lemma}

\begin{proof}
    We use a similar sieve 
    \begin{equation}
        \mathcal{F}_n=\left\{\Omega\in \mathcal{B}^\Omega_n: ||\Omega||_1 \le 8C_3q \right\}
    \end{equation}
    The Sieve covering can be checked with results in~\eqref{eqn:sieve_Omega_bound_mssl}. 
\end{proof}

\begin{lemma}[Test condition]
\label{lemma:test-condition}
    We denote the dimension recover event as  $\mathcal{B}^\Omega_n:=\{\Omega\succ \tau I:|\nu_{\omega}(\Omega)|<C_3'\max\{q,s_0^\Omega\}\}$. For the sieve $\mathcal{F}_n:=\left\{\Omega\in \mathcal{B}^\Omega_n: ||\Omega||_1 \le 8C_3q \right\}$ there exist some constants $C_2>C_1+2$:
\begin{equation}
    \Pi(\mathcal{F}_n^c)\le \exp(-C_2n\epsilon_n^{'2}).
\end{equation}
There exists tests $\varphi_n$, such that for some $M_2>C_1+1$:
\begin{equation}
    \label{eqn:test_condition_gssl}
    \begin{aligned}
        &\mathbb{E}_{f_0}\varphi_n\lesssim e^{-M_2n\epsilon^2/2}\\
        \sup_{f\in \mathcal{F}_n:\rho(f_0,f)>M_2n\epsilon_n^2}&\mathbb{E}_{f}(1-\varphi_n)\lesssim e^{-M_2n\epsilon_n^2}
    \end{aligned}
\end{equation}
where $f=\prod_{i=1}^n \mathcal{N}(0,\Omega^{-1})$ while $f_0=\prod_{i=1}^n \mathcal{N}(0,\Omega_0^{-1})$
\end{lemma}

\begin{proof}
    Construct Neyman-Pearson test on
    $|||\Omega|||_2\le ||\Omega||_1\le 8C_3q$,
         and 
        $|||\Omega_1-\Omega|||_2\le ||\Omega_1-\Omega||_1\le \frac{1}{8C_3nq^{3/2}}. $, the Neyman-Pearson test has vanishing type-II error rate, same as in mSSL, then we take supremim over these sets and using the same covering argument as in~\eqref{eq:covering_Omega}. 
\end{proof}

From log-affinity we can have parameter recovery similar to mSSL.

\begin{theorem}[Contraction of gSSL]
    For Gaussian graphical model $Y_i\sim \N(0,\Omega^{-1})$, with prior on $\Omega$ following the SSL prior and Assumptions A1, A5, and 
    \begin{itemize}
        \item[A3']: $\log(n)\lesssim \log(q)$ and $\max\{q,s_0^\Omega\}\log q/n\to 0$
    \end{itemize}
    we have 
    \begin{align}
        \sup_{\Omega\in\mathcal{H}_0}\E_0 \Pi\left(\Omega:||\Omega-\Omega_0||_F^2\ge M_1\epsilon_n^{'2}|Y_1,\dots,Y_n\right) &\longrightarrow 0
    \end{align}
    where $\epsilon'_n = \max\{q,s_0^\Omega\}\log q/n$
\end{theorem}

\begin{proof}
    Use Theorem~\ref{thm:log-affinity-gssl} and same bounds on $\Omega$ as Equation (14) in the main text.
\end{proof}

\section{Details and additional results in experiments}
\label{appendix:experiments}
\subsection{Details on covariance structures}
\label{appendix:cov_structure}
For each choice of $(n,p,q),$ we considered five different $\Omega$'s: (i) an AR(1) model for $\Omega^{-1}$ so that $\Omega$ is tri-diagonal; (ii) an AR(2) model for $\Omega^{-1}$ so that $\omega_{k,k'} = 0$ whenever $\lvert k - k'\vert > 2$; (iii) a block model in which $\Omega$ is block-diagonal with two dense $q/2 \times q/2$ diagonal blocks; (iv) a small-world network; and (v) a tree grap.

In the AR(1) model we set $(\Omega^{-1})_{k,k'} = 0.7^{\lvert k - k' \rvert}$ so that $\omega_{k,k'} = 0$ whenever $\lvert k - k' \vert > 1.$
In the AR(2) model, we set $\omega_{k,k} = 1, \omega_{k-1,k} = \omega_{k,k-1} = 0.5,$ and $\omega_{k-2,k} = \omega_{k,k-2} = 0.25.$
For the block model, we partitioned $\Sigma = \Omega^{-1}$ into 4 $q/2 \times q/2$ blocks and set all entries in the off-diagonal blocks of $\Sigma$ to zero. 
We then set $\sigma_{k,k} = 1$ and $\sigma_{k,k'} = 0.5$ for $1 \leq k \neq k' \leq q/2$ and for $q/2 + 1 \leq k \neq k' \leq q.$
For the small-world and tree networks, we drew $\Omega$ from a G-Wishart distribution \citep{Roverato2002, Lenkoski2013} with three degrees of freedom and identity scale matrix after drawing the underlying graphical structure.
We drew the small-world graph using the Watts-Strogatz \citep{WattsStrogatz1998} model with a single community and rewiring probability of 0.1.
And we drew the tree by running a loop-erased random walk on the complete graph.
These settings are identical to those used in \cite{shen2022cgSSL}.

\subsection{Other dimensions}

Here we provide experiments on different dimensions with $(n,p,q)=(100, 10, 10)$ and $(n,p,q)=(100, 20, 30)$.

\begin{table}[H]
    \centering
        \caption{Support recovery and uncertainty interval coverage for $B$ when $(n,p,q) = (100, 20, 30)$ average across 100 simulated datasets. Best performance is bold-faced.}
    \label{tab:simulation_on_debias_b_small}
    \scriptsize
    \begin{tabular}{lccccccc}
        \hline
        Method & SEN & SPEC & PREC & F1 & Relative length& Coverage & Interval score  \\ \hline
        \multicolumn{8}{c}{$AR(1)$ model} \\ \hline
        \texttt{debiasing} & \textbf{0.86} & 0.95 & 0.81 & 0.83 & 1 & 0.95 & 7.7 \\
        \texttt{WLB} & 0.78 & \textbf{1} & \textbf{1} & 0.87 & \textbf{0.11} & 0.95 & \textbf{1.1} \\
        \texttt{WBB} & 0.78 & \textbf{1} & \textbf{1} & 0.88 & 0.35 & \textbf{0.98} & 3.9 \\
        \texttt{BB-jitter} & 0.79 & \textbf{1} & \textbf{1} & 0.88 & 0.52 & \textbf{0.98} & 3 \\
        \texttt{mSSL} & 0.82 & \textbf{1} & \textbf{1} & \textbf{0.90} & NA & NA & NA \\
        \hline
        \multicolumn{8}{c}{$AR(2)$ model} \\ \hline
        \texttt{debiasing} & \textbf{0.82} & 0.95 & 0.8 & 0.81 & 1 & 0.95 & 10 \\
        \texttt{WLB} & 0.59 & \textbf{1} & \textbf{1} & 0.74 & \textbf{0.16} & 0.93 & \textbf{2.3} \\
        \texttt{WBB} & 0.60 &\textbf{1} & \textbf{1} & 0.75 & 0.36 & 0.96 & 5.7 \\
        \texttt{BB-jitter} & 0.63 & \textbf{1} & \textbf{1} & 0.77 & 0.56 & \textbf{0.97} & 4.5\\
        \texttt{mSSL} & 0.72 & \textbf{1} & \textbf{1} & \textbf{0.84} & NA & NA & NA \\
        \hline
        \multicolumn{8}{c}{Block model} \\ \hline
        \texttt{debiasing} &\textbf{0.86} & 0.95 & 0.82 & 0.84 & 1 & 0.95 & 7.7 \\
        \texttt{WLB} & 0.75 & \textbf{1} & \textbf{1} & 0.86 & \textbf{0.13} & 0.94 & \textbf{1.3} \\
        \texttt{WBB} & 0.76 & \textbf{1} & \textbf{1} & 0.86 & 0.36 & 0.97 & 4.1 \\
        \texttt{BB-jitter} & 0.77 & \textbf{1} & \textbf{1} & 0.87 & 0.34 & \textbf{0.98} & 3.2 \\
        \texttt{mSSL} & 0.81 & \textbf{1} & 0.99 & \textbf{0.89} & NA & NA & NA \\
        \hline
        \multicolumn{8}{c}{Small world} \\ \hline
        \texttt{debiasing} & \textbf{0.88} & 0.95 & 0.81 & 0.85 & 1 & 0.95 & 6.4 \\
        \texttt{WLB} &0.78 & \textbf{1} & \textbf{1} & 0.88 & \textbf{0.14} & 0.95 & \textbf{1.0}\\
        \texttt{WBB} &0.78 & \textbf{1} & \textbf{1} & 0.88 & 0.34 & \textbf{0.98} & 3.6 \\
        \texttt{BB-jitter} &0.80 & \textbf{1} & \textbf{1} & 0.89 & 0.45 & \textbf{0.98} & 2.9\\
        \texttt{mSSL} &  0.83 & \textbf{1} & \textbf{1} & \textbf{0.90}  & NA & NA & NA \\
        \hline
        \multicolumn{8}{c}{Tree} \\ \hline
        \texttt{debiasing} & \textbf{0.90} & 0.95 & 0.82 & 0.86 & 1 & 0.95 & 5.3 \\
        \texttt{WLB} &0.78 & \textbf{1} & \textbf{1} & 0.87 & \textbf{0.17} & 0.96 & \textbf{1.0} \\
        \texttt{WBB} &0.78 & \textbf{1} & \textbf{1} & 0.88 & 0.44 & \textbf{0.98} & 3.4 \\
        \texttt{BB-jitter} &0.80 & \textbf{1} & \textbf{1} & 0.89 & 0.55 & \textbf{0.98} & 2.7\\
        \texttt{mSSL} & 0.84 & \textbf{1} & \textbf{1} & \textbf{0.91} & NA & NA & NA \\
        \hline
    \end{tabular}
\end{table}

\begin{table}[H]
    \centering
        \caption{Support recovery and uncertainty interval coverage for $\Omega$ when $(n,p,q) = (100, 20, 30)$ average across 100 simulated datasets. Best performance is bold-faced.}    \label{tab:simulation_on_debias_omega_small}
    \scriptsize
    \begin{tabular}{lccccccc}
        \hline
        Method & SEN & SPEC & PREC & F1 & Length& Coverage & Interval score  \\ \hline
        \multicolumn{8}{c}{$AR(1)$ model} \\ \hline
        \texttt{debiasing} & \textbf{1} & 0.96 & 0.64 & 0.78 & 1 & 0.96 & 23 \\
        \texttt{WLB} & \textbf{1} & \textbf{1} & \textbf{1} & \textbf{1} & \textbf{0.052} & \textbf{1} & \textbf{1.2} \\
        \texttt{WBB} & \textbf{1} & \textbf{1} &\textbf{1} & \textbf{1}& 0.16 & \textbf{1} & 4.6 \\
        \texttt{BB-jitter} & \textbf{1} & \textbf{1} & \textbf{1} & \textbf{1} & 0.38 & \textbf{1} & 3.6 \\
        \texttt{mSSL} & \textbf{1} & \textbf{1} & 0.99 & 1 & NA & NA & NA\\
        \hline
        \multicolumn{8}{c}{$AR(2)$ model} \\ \hline
        \texttt{debiasing} & \textbf{0.66} & 0.92 & 0.56 & \textbf{0.6} & 1 & 0.82 & 6.0\\
        \texttt{WLB} & 0.055 & \textbf{1} & 0.31 & 0.085 & \textbf{0.017} & 0.87 & \textbf{1.9} \\
        \texttt{WBB} & 0.055 & \textbf{1} & 0.32 & 0.085 & 0.73 & \textbf{0.97} & 5.0 \\
        \texttt{BB-jitter} & 0.057 & \textbf{1} & 0.39 & 0.089 & 0.26 & \textbf{0.97} & 3.7\\
        \texttt{mSSL} & 0.068 & \textbf{1} & \textbf{0.78} & 0.11 & NA & NA & NA \\
        \hline
        \multicolumn{8}{c}{Block model} \\ \hline
        \texttt{debiasing} & 0.17 & 0.93 & 0.65 & 0.25 & 1 & \textbf{0.90} & 17 \\
        \texttt{WLB} & 0.17 & \textbf{1} & 0.94 & 0.29 & \textbf{0.065} & 0.52 & \textbf{5.4} \\
        \texttt{WBB} & 0.14 & \textbf{1} & \textbf{0.95} & 0.24 & 0.32 & 0.81 & 6.5 \\
        \texttt{BB-jitter} & 0.17 & \textbf{1} & \textbf{0.95} & 0.29 & 1.2 & 0.85 & 6.4 \\
        \texttt{mSSL} & \textbf{0.18} & \textbf{1} & \textbf{0.95} & \textbf{0.30} & NA & NA & NA\\
        \hline
        \multicolumn{8}{c}{Small world} \\ \hline
        \texttt{debiasing} & \textbf{0.64} & 0.96 & 0.54 & 0.58 & 1 & 0.96 & 29 \\
        \texttt{WLB} &0.51 & \textbf{1} & 0.94 & \textbf{0.66} & \textbf{0.034} & 0.92 & \textbf{2.1}\\
        \texttt{WBB} &0.49 & \textbf{1} & \textbf{0.96} & 0.65 & 0.4 & 0.94 & 5.7 \\
        \texttt{BB-jitter} &0.51 & \textbf{1} & 0.95 & \textbf{0.66} & 0.41 & \textbf{0.97} & 5.4\\
        \texttt{mSSL} &  0.51 & \textbf{1}& 0.94 & \textbf{0.66} & NA & NA & NA \\
        \hline
        \multicolumn{8}{c}{Tree} \\ \hline
        \texttt{debiasing} & \textbf{0.71} & 0.95 & 0.53 & \textbf{0.61} & 1 & 0.95 & 33 \\
        \texttt{WLB} &0.45 & \textbf{1} & 0.93 & 0.60 & \textbf{0.033} & 0.89 & \textbf{2.9} \\
        \texttt{WBB} &0.43 & \textbf{1} & \textbf{0.96} & 0.59 & 0.18 & 0.91 & 5.6\\
        \texttt{BB-jitter} &0.45 & \textbf{1} & 0.93 & 0.60 & 0.21 & \textbf{0.97} & 5.7\\
        \texttt{mSSL} & 0.45 & \textbf{1} & 0.93 & \textbf{0.61} & NA & NA & NA \\
        \hline
    \end{tabular}
\end{table}


\begin{table}[H]
    \centering
        \caption{Support recovery and uncertainty interval coverage for $B$ when $(n,p,q) = (100, 10, 10)$ average across 100 simulated datasets. Best performance is bold-faced.}
        \label{tab:simulation_on_debias_b_tiny}
    \scriptsize
    \begin{tabular}{lccccccc}
        \hline
        Method & SEN & SPEC & PREC & F1 & Relative length& Coveage & Interval score  \\ \hline
        \multicolumn{8}{c}{$AR(1)$ model} \\ \hline
        \texttt{debiasing} & \textbf{0.97} & 0.95 & 0.83 & 0.89 & 1 & 0.95 & 7.6 \\
        \texttt{WLB} & 0.88 & \textbf{1} & \textbf{1} & 0.93 &\textbf{0.13} & 0.97 & \textbf{1.2}  \\
        \texttt{WBB} & 0.88 & \textbf{1} & \textbf{1} & 0.94 & 0.37 & \textbf{0.98} & 4.0 \\
        \texttt{BB-jitter} & 0.92 & \textbf{1} & \textbf{1} & 0.96 & 0.6 & \textbf{0.98} & 3.2 \\
        \texttt{mSSL} & 0.95 & \textbf{1} & \textbf{1} & \textbf{0.98} & NA & NA & NA \\
        \hline
        \multicolumn{8}{c}{$AR(2)$ model} \\ \hline
        \texttt{debiasing} & \textbf{0.93} & 0.95 & 0.82 & 0.87 & 1 & 0.95 & 9.9 \\
        \texttt{WLB} & 0.61 & \textbf{1} & \textbf{1} & 0.75 & \textbf{0.17} & 0.95 & \textbf{2.3} \\
        \texttt{WBB} & 0.61 & \textbf{1} & \textbf{1} & 0.75 & 0.68 & \textbf{0.97} & 5.6 \\
        \texttt{BB-jitter} & 0.65 & \textbf{1} & \textbf{1} & 0.79 & 0.57 & \textbf{0.97} & 4.4\\
        \texttt{mSSL} & 0.82 & \textbf{1} & \textbf{1} & \textbf{0.90} & NA & NA & NA \\
        \hline
        \multicolumn{8}{c}{Block model} \\ \hline
        \texttt{debiasing} &\textbf{0.97} & 0.95 & 0.84 & 0.90 & 1 & 0.95 & 7.5 \\
        \texttt{WLB} & 0.81 & \textbf{1} & \textbf{1} & 0.89 & \textbf{0.16} & 0.97 & \textbf{1.5} \\
        \texttt{WBB} & 0.82 & \textbf{1} & \textbf{1} & 0.90 & 0.71 & \textbf{0.98} & 4.4\\
        \texttt{BB-jitter} & 0.86 & \textbf{1} & \textbf{1} & 0.92 & 0.38 & \textbf{0.98} & 3.5 \\
        \texttt{mSSL} & 0.92 & 1 & 0.99 & \textbf{0.96} & NA & NA & NA\\
        \hline
        \multicolumn{8}{c}{Small world} \\ \hline
        \texttt{debiasing} & \textbf{0.98} & 0.95 & 0.83 & 0.90 & 1 & 0.95 & 5.6 \\
        \texttt{WLB} &0.87 & \textbf{1} & \textbf{1} & 0.93 & \textbf{0.19} & 0.97 & \textbf{1.2}\\
        \texttt{WBB} &0.88 & \textbf{1} & \textbf{1} & 0.94 & 0.28 & \textbf{0.98} & 3.4 \\
        \texttt{BB-jitter} &0.9 & \textbf{1} & \textbf{1} & 0.95 & 0.48 & \textbf{0.98} & 2.8\\
        \texttt{mSSL} &  0.95 & \textbf{1} & \textbf{1} & \textbf{0.97}  & NA & NA & NA \\
        \hline
        \multicolumn{8}{c}{Tree} \\ \hline
        \texttt{debiasing} & \textbf{1} & 0.95 & 0.83 & 0.91 & 1 & 0.95 & 4.3 \\
        \texttt{WLB} &0.94 & \textbf{1} & \textbf{1} & 0.97 & \textbf{0.20} & 0.98 & \textbf{0.94} \\
        \texttt{WBB} &0.94 & \textbf{1} & \textbf{1} & 0.97 & 0.27 & \textbf{0.99} & 3.0 \\
        \texttt{BB-jitter} &0.97 & \textbf{1} & \textbf{1} & \textbf{0.99} & 0.60 & \textbf{0.99} & 2.5\\
        \texttt{mSSL} & 0.98 & \textbf{1} & \textbf{1} & \textbf{0.99}  & NA & NA & NA \\
        \hline
    \end{tabular}
\end{table}

\begin{table}[H]
    \centering
        \caption{Support recovery and uncertainty interval coverage for $\Omega$ when $(n,p,q) = (100, 10, 10)$ average across 100 simulated datasets. Best performance is bold-faced.}
    \label{tab:simulation_on_debias_omega_tiny}
    \scriptsize
    \begin{tabular}{lccccccc}
        \hline
        Method & SEN & SPEC & PREC & F1 & Length& Coverage & Interval score  \\ \hline
        \multicolumn{8}{c}{$AR(1)$ model} \\ \hline
        \texttt{debiasing} & \textbf{1} & 0.95 & 0.85 & 0.92 & 1 & 0.95 & 22 \\
        \texttt{WLB} & \textbf{1} & \textbf{1} & \textbf{1} & \textbf{1} & \textbf{0.15} & 0.99 & \textbf{3.4} \\
        \texttt{WBB} & \textbf{1} & \textbf{1} & \textbf{1} & \textbf{1}& 0.76 & 0.99 & 6.9 \\
        \texttt{BB-jitter} & \textbf{1} & \textbf{1} & \textbf{1} & \textbf{1} & 0.51 & \textbf{0.99} & 5.7\\
        \texttt{mSSL} & \textbf{1} & \textbf{1} & \textbf{1} & \textbf{1} & NA & NA & NA\\
        \hline
        \multicolumn{8}{c}{$AR(2)$ model} \\ \hline
        \texttt{debiasing} & \textbf{0.77} & 0.91 & \textbf{0.85} & \textbf{0.81} & 1 & 0.81 & 8.1\\
        \texttt{WLB} & 0.45 & 0.99 & 0.72 & 0.54 & \textbf{0.17} & 0.75 & \textbf{4.6} \\
        \texttt{WBB} & 0.41 & \textbf{1} & 0.74 & 0.52 & 0.36 & 0.87 & 6.2 \\
        \texttt{BB-jitter} & 0.45 & 0.99 & 0.73 & 0.54 & 0.29 & \textbf{0.88} & 5.2\\
        \texttt{mSSL} & 0.47 & 0.99 & 0.76 & 0.56 & NA & NA & NA \\
        \hline
        \multicolumn{8}{c}{Block model} \\ \hline
        \texttt{debiasing} & 0.48 & 0.94 & 0.89 & 0.61 & 1 & \textbf{0.94} & 14 \\
        \texttt{WLB} &0.51 &\textbf{1} & \textbf{1} & \textbf{0.68} & 0.20 & 0.65 & \textbf{7.4} \\
        \texttt{WBB} & 0.44 & \textbf{1} & \textbf{1} & 0.61 & 0.25 & 0.83 & 8.8 \\
        \texttt{BB-jitter} & 0.50 & \textbf{1} & \textbf{1} & 0.66 & 1.0 & 0.85 & 8.3 \\
        \texttt{mSSL} & \textbf{0.52} & \textbf{1} & \textbf{1} & \textbf{0.68} & NA & NA & NA\\
        \hline
        \multicolumn{8}{c}{Small world} \\ \hline
        \texttt{debiasing} & \textbf{0.76} & 0.95 & 0.84 & \textbf{0.79} & 1 & \textbf{0.95} & 31 \\
        \texttt{WLB} &0.51 & \textbf{1} & \textbf{0.99} & 0.67 & \textbf{0.097} & 0.88 & \textbf{6.8} \\
        \texttt{WBB} &0.49 & \textbf{1} & \textbf{0.99} & 0.65 & 0.68 & 0.91 & 10 \\
        \texttt{BB-jitter} &0.51 & \textbf{1} & \textbf{0.99} & 0.67 & 0.57 & 0.91 & 9.7\\
        \texttt{mSSL} &  0.52 & \textbf{1} & \textbf{0.99} & 0.67  & NA & NA & NA \\
        \hline
        \multicolumn{8}{c}{Tree} \\ \hline
        \texttt{debiasing} & \textbf{0.55} & 0.95 & 0.78 & \textbf{0.64} & 1 & 0.95 & 35 \\
        \texttt{WLB} &0.33 & \textbf{1} & \textbf{0.97} & 0.49 & \textbf{0.063} & 0.8 & \textbf{6.2}\\
        \texttt{WBB} &0.33 & \textbf{1} & \textbf{0.97} & 0.48 & 0.53 & 0.84 & 8.9 \\
        \texttt{BB-jitter} &0.33 & \textbf{1} & \textbf{0.97} & 0.49 & 0.24 & 0.90 & 9.2\\
        \texttt{mSSL} &  0.33 & \textbf{1} & \textbf{0.97} & 0.49  & NA & NA & NA \\
        \hline
    \end{tabular}
\end{table}

\end{document}